\documentclass[12pt,leqno]{article}

\usepackage{amsthm}%,showidx}

\usepackage{amsmath}
\usepackage{amscd}
\usepackage{amsmath}
\usepackage{amssymb}
\usepackage{latexsym}
\usepackage[title]{appendix}
\usepackage{color}

\makeatletter

\@addtoreset{equation}{section}
\makeatother

\newtheorem{thm}{Theorem}[section]
\newtheorem{pr}[thm]{Proposition}
\newtheorem{df}[thm]{Definition}
\newtheorem{lm}[thm]{Lemma}
\newtheorem{cor}[thm]{Corollary}

\newtheorem{rmk}[thm]{Remark}
\newtheorem{ex}[thm]{Example}
\newtheorem{df-lm}[thm]{Definition-Lemma}
\newtheorem{thmA}{Theorem}[section]
\newtheorem{prA}[thmA]{Proposition}
\newtheorem{dfA}[thmA]{Definition}
\newtheorem{lmA}[thmA]{Lemma}
\newtheorem{exA}[thmA]{Example}

\newcommand{\C}{\mathbb{C}}
\newcommand{\rk}{{\rm rk}}

\newcommand{\F}{\mathbb{F}}

\newcommand{\Ql}{\overline{\mathbb{Q}_\ell}}
\newcommand{\G}{\mathbb{G}}
\newcommand{\Qp}{\mathbb{Q}_p}
\newcommand{\Zp}{\mathbb{Z}_p}

\newcommand{\Z}{\mathbb{Z}}

\input xy \xyoption{all} \CompileMatrices

\begin{document}

\title{Quadratic $\ell$-adic sheaf and its Heisenberg group}
\author{Daichi Takeuchi
\thanks{RIKEN, Center for Advanced Intelligence Project AIP, Mathematical Science Team
\texttt{Email: daichi.takeuchi@riken.jp}}}
\date{}
\maketitle

\begin{abstract}
In this paper, we introduce a new class of $\ell$-adic sheaves, which we call quadratic $\ell$-adic sheaves, on connected unipotent commutative algebraic groups over  finite fields. They are sheaf-theoretic enhancements of quadratic forms on finite abelian groups in the spirit of the function-sheaf dictionary. We show that a certain finite Heisenberg group acts on a quadratic sheaf and that the cohomology of the quadratic sheaf gives an irreducible representation of the group. 

We also compute the Frobenius eigenvalues of the cohomology groups. As a byproduct, we find a large number of examples of affine supersingular varieties. 
\end{abstract}
\section{Introduction}
Classically, quadratic Gauss sums are defined as follows. Let $\psi\colon\F_q\to\C^\times$ be a non-trivial character on a finite field $\F_q$ with $q$ elements. Then the associated quadratic Gauss sum is defined to be $\sum_{x\in\F_q}\psi(x^2)$. This kind of sums naturally appears in many examples, such as in the $L$-functions of some classes of algebraic varieties over $\F_q$. It also appears when one computes the local epsilon factors of characters that are wildly ramified. However, 
when we consider this kind of Gauss sums, we always make the assumption  that $q$ is odd: when $q$ is even, the sum given by the same formula is reduced to $0$ as the map $x\mapsto \psi(x^2)$ is again a non-trivial character. This sometimes forces us to make division into cases according to the parity of $q$ when we do a concrete computation. 
On the other hand, it is known that there is a generalization of the formulation of quadratic Gauss sum in which one can also treat the case when $q$ is even. This is done by replacing $\F_q$ with a finite abelian group of arbitrary order and $\psi(x^2)$ with a function with a certain property, which is called as a quadratic form valued in $\mathbb{C}^\times$: we review this in the appendix A of this paper. This reformulation can be useful to express quantities in a uniform way, although the division into cases {\it is}  necessary if one tries to compute them more specifically. 

Let $\ell$ be a prime number which is invertible in finite fields that appear in what follows. The aim of this paper is to introduce a new class of $\ell$-adic sheaves, which we call as {\it quadratic $\ell$-adic sheaves}, which can be viewed as sheaf-theoretic enhancements of quadratic forms on finite abelian groups, in the spirit of the function-sheaf dictionary of $\ell$-adic sheaves. 

We describe the contents of this paper more precisely. We fix an algebraic closure $\Ql$ of the $\ell$-adic field $\mathbb{Q}_\ell$. By an $\ell$-adic sheaf, we mean a $\Ql$-sheaf. Let $U$ be a connected unipotent commutative algebraic group over $\F_q$ (where the parity of $q$ now can be either of odd or even). A quadratic sheaf is defined as an invertible $\ell$-adic sheaf ${\cal Q}$ on $U$ satisfying a certain condition that guarantees that the map $U(\F_q)\to\Ql^\times$ obtained by taking the Frobenius traces is a quadratic form (in the sense of Definition \ref{quadcl} if one replaces $\mathbb{C}$ with $\Ql$) on a finite abelian group $U(\F_q)$. For the precise definition of quadratic sheaves, see Definition \ref{quaddef}.1. 

We actually consider the perfection $U^{\rm perf}$ of $U$, rather than $U$ itself, as the base scheme on which quadratic sheaves live.  This is because the study  on quadratic sheaves done in this paper heavily relies on the duality theorem of perfect unipotent commutative groups \cite{Beg}; as the \'etale topos of $U$ is canonically equivalent to that of $U^{\rm perf}$, all results on \'etale sheaves on $U^{\rm perf}$ can be automatically rephrased into those on $U$. However, staying in the category of perfect schemes makes statements much simpler. For this reason, in the main body of this paper, we use the symbol $U$ for perfect unipotent groups rather than algebraic ones.

Our first main result is a computation of the ranks of the cohomology groups of $\cal Q$. 
\begin{thm}(Theorem \ref{HD})\label{introHD}
Suppose that $\cal Q$ is a quadratic $\Ql$-sheaf on $U$ that is isogeneous of degree $p^{2r}$ (Definition \ref{quaddef}.2). Let $d$ be the dimension of $U$. Then the cohomology group $H^i_c(U_{\overline{\F_q}},{\cal Q})$ vanishes unless $i=d$. Moreover, the rank of $H^d_c(U_{\overline{\F_q}},{\cal Q})$ is equal to $p^r$.  
\end{thm}

This computation leads us to the following version of Hasse--Davenport relation on Gauss sums. 
\begin{cor}(Corollary \ref{genHD} for a more general statement)
Let $\cal Q$ be a quadratic $\Ql$-sheaf on $U$ that is non-degenerate (i.e., isogeneous of degree $1$). For an integer $n\geq1$, let $t_{{\cal Q},q^n}$ denote the function $U(\F_{q^n})\to\Ql^\times$ given by the Frobenius traces. Then we have 
\[
(-1)^d\sum_{x\in U(\F_{q^n})}t_{{\cal Q},q^n}(x)=((-1)^d
\sum_{x\in U(\F_{q})}t_{{\cal Q},q}(x))^n
\]
where $d=\dim U$. 
\end{cor}
Along the way of investigating quadratic sheaves, we find that a finite Heisenberg group $H$ (i.e., a two-step nilpotent group) can be naturally attached to an isogeneous quadratic sheaf. Very briefly, this group is constructed as follows. Let $U^\ast$ be the dual group of $U$ (defined as a perfect group scheme). Then the bimultiplicative sheaf associated to the quadratic sheaf $\cal Q$ gives a group homomorphism $U\to U^\ast$. The condition that $\cal Q$ is isogeneous means that this morphism is an isogeny. Let $C\to U$ be the Galois covering that trivializes $\cal Q$. Then $H$ is defined to be the Galois group of the composite map $C\to U\to U^\ast$. See the subsection \ref{H} for a detailed account. 

The cohomology $H^d_c(U_{\overline{\F_q}},{\cal Q})$ is also interesting from the viewpoint of the representation theory for $H$. 
\begin{thm}(Theorem \ref{irrH})
Let $\cal Q$ be an isogeneous quadratic $\Ql$-sheaf on $U$. Let $H$ be the Heisenberg group constructed in the subsection \ref{H} of the main body of this paper. Then $H^d_c(U_{\overline{\F_q}},{\cal Q})$ is irreducible  as an $H$-representation. 
\end{thm}

In Section \ref{secH}, we study in detail this Heisenberg group, especially its geometric action on algebraic varieties, and deduce several consequences. Among them, we prove the following result. 
\begin{thm}(Corollary \ref{nondegdes}) 
Let $\cal Q$ be an isogeneous quadratic $\Ql$-sheaf on $U$. Then, after replacing the base field $\F_q$ by a finite extension if necessary, one can find another connected commutative unipotent algebraic group $V$ with an abelian \'etale isogeny $f\colon U\to V$ and a non-degenerate quadratic $\Ql$-sheaf ${\cal Q}'$ on $V$ such that the pullback $f^\ast{\cal Q}'$ is isomorphic to $\cal Q$. 
\end{thm}

This result is useful to study the Frobenius action on 
$H^d_c(U_{\overline{\F_q}},{\cal Q})$. Let $K$ denote the kernel of $f$. 
Up to a finite base field extension, this theorem gives us  the following direct sum decomposition 
\[H^d_c(U_{\overline{\F_q}},{\cal Q})\cong \bigoplus_{\psi\in K^\vee}H^d_c(V_{\overline{\F_q}},{\cal Q}'\otimes{\cal L}_\psi)
\]
where $K^\vee$ denotes the character group ${\rm Hom}(K(\overline{\F_q}),\Ql^\times)$ and ${\cal L}_\psi$ is the invertible $\Ql$-sheaf on $V$ given by $\psi$. Since ${\cal Q}'\otimes{\cal L}_\psi$ are non-degenerate quadratic sheaves, a computation of $H^d_c(U_{\overline{\F_q}},{\cal Q})$ can be reduced to the non-degenerate case. In the non-degenerate case, the rank of the cohomology is equal to $1$ according to Theorem \ref{introHD}. Then Grothendieck's trace formula and an evaluation of Gauss sum (Proposition \ref{GS}) imply that any eigenvalue of the $q$-th power Frobenius ${\rm Fr}_q$ on $H^d_c(U_{\overline{\F_q}},{\cal Q})$ is equal to $\sqrt{q^d}$ up to a root of unity. 

From this computation, one can show that the covering space $C\to U$ that trivializes a quadratic sheaf is 
{\it supersingular}; here we say that an affine smooth variety $X$ over $\F_q$ is supersingular if all the eigenvalues of ${\rm Fr}_q$ on $H^i_c(X_{\overline{\F_q}},\Ql)$ are equal to $\sqrt{q^i}$ up to roots of unity for any $i$. Using quadratic sheaves, one can produce a large number of examples of affine supersingular varieties. 
See the subsection \ref{AffSS} for a detailed discussion. In the appendix \ref{exampleC}, we compute some examples of such varieties in low dimensional cases. In the curve case, the Van der Geer--Van der Vlugt curves appear. This family of curves is extensively studied by many authors (for example, \cite{GV}, \cite{BHetal}, \cite{TT}). Thus the theory of quadratic sheaves can be a useful tool to study such curves. 

As explained above, an isogeneous quadratic sheaf $\cal Q$ has an action by  a Heisenberg group. More interestingly, one can construct another group which acts on $\cal Q$: thus it turns out that the cohomology group 
$H^d_c(U_{\overline{\F_q}},{\cal Q})$ is the underlying
space of a representation of a large group. The appendix \ref{unitary} is devoted to explaining this observation. Although the representations so obtained remain mysterious in this paper, it 
 seems that they are important for the following reason. 
It is observed by Imai--Tsushima \cite{IT} and Tsushima \cite{Tsu} that the Heisenberg--Weil representations of some finite reductive groups can be realized geometrically. Their constructions can be regarded as special cases of ours (Example \ref{HWex}). It would be an interesting problem to study this group and the representation in detail. 

We explain the construction of this paper. As mentioned above, we mainly use  perfect unipotent commutative group schemes (Definition \ref{punip}) rather than algebraic ones. 
In Section $2$, we review definitions and basic results on them. We also prove some preliminary results in this section. In Section $3$, we define the quadratic sheaf by using the terminologies recalled in the previous section. We investigate quadratic sheaves in detail in the following sections. In the appendix A, we record results on quadratic forms on finite abelian groups, among which Proposition \ref{GS} is important for us. 
Although this result seems well-known, the author cannot find an appropriate reference. Thus we give a complete proof of this proposition for readers' convenience. This proof can be  a guideline for the proof of Theorem \ref{HD}.

\section*{Convention and Notation}
We collect convention and notation that we use throughout this paper. 
Let $p$ be a prime number and 
$k$ be a perfect field of characteristic $p$. We write 
$\Lambda$ for the coefficient ring of \'etale sheaves.  We  assume that it is either of $\overline{\mathbb{Q}_\ell}$ for a prime number $\ell\neq p$ or a noetherian strictly local ring in which some integer $\geq1$ invertible in $k$ is zero (e.g. $W(\overline{\F_\ell})/(\ell^n)$); our arguments can be applied to the both cases without any modification. 

\begin{itemize}
\item For a scheme $S$, we write $Sch/S$ for the category of $S$-schemes. When $S$ is a perfect $\F_p$-scheme, we write $PSch/S$ for the full subcategory consisting of perfect $S$-schemes. Presheaves mean set-valued contravariant functors. 
\item Let $S$ be a perfect $\F_p$-scheme. For $X\in Sch/S$, we write $X^{\rm perf}\in PSch/S$ for its perfection, i.e., $\varprojlim_{\rm Frob}X$. For a presheaf $F\colon (Sch/S)^{\rm op}\to ({\rm Sets})$, we write $F^{\rm perf}$ for the restriction $F|_{PSch/S}$ and also call it the perfection. 
If we identify $X\in Sch/S$ with the presheaf on $Sch/S$ given by the Yoneda embedding, then the perfection as a presheaf is represented by the perfect scheme $X^{\rm perf}$; the notation is compatible with the Yoneda embedding. 
\item For an $S$-scheme $X$ and a morphism of schemes $T\to S$, we write $X_T$ for the $T$-scheme $X\times_ST$. 
\item For a morphism $f\colon X\to Y$ of schemes and an \'etale sheaf $\cal L$ on $Y$, we write ${\cal L}(f)$ for $f^\ast{\cal L}$. 
\item For a finite abelian group $M$ such that $|M|$ is invertible in $\Lambda$, we write $M^\vee$ for the character group ${\rm Hom}(M,\Lambda^\times)$. 
\end{itemize}

\tableofcontents

\section{Preliminaries}
In this section, we collect definitions and results which are necessary to define and investigate quadratic sheaves.

\subsection{Perfect unipotent commutative groups}
Let $k$ be a perfect field of characteristic $p>0$. 
\begin{df}\label{punip}
A group $k$-scheme  is said to be {\rm perfect unipotent commutative} if there exists a connected affine smooth unipotent commutative group $k$-scheme $U_0$ such that it is isomorphic to the perfection $U_0^{\rm perf}$. 
\end{df}
Let $U$ be a perfect unipotent group $k$-scheme. By an algebraic model of $U$, we mean 
an affine smooth unipotent group $k$-scheme $U_0$ such that $U_0^{\rm perf}\cong U$. Note that the underlying space of $U$ as well as the \'etale topos of $U$ can be identified with those of $U_0$, as perfection does not change the underlying topological space and the \'etale site. 

The following lemma is useful to reduce computations for perfect unipotent groups to those for algebraic models. 
\begin{lm}\label{redalg}
Let $f\colon U\to V$ be a morphism of perfect unipotent commutative group schemes over $k$. Then there exist algebraic models $U_0,V_0$ of $U,V$ respectively together with a $k$-morphism $f_0\colon U_0\to V_0$ of group $k$-schemes such that the perfection of $f_0$ is equal to $f$. 
\end{lm}
In the sequel, we call such a morphism $f_0$ as {\it an algebraic model} of $f$. 
\proof{
For an affine $k$-scheme $X$, let ${\cal O}(X)$ denote its affine ring. 

Let $U_0,V_0$ be unipotent affine smooth group schemes over $k$ such that their perfections are isomorphic to $U,V$ respectively. We have the following diagram
 \begin{equation*}
 {\cal O}(V_0)\hookrightarrow{\cal O}(V)\xrightarrow{f^\ast}{\cal O}(U)\hookleftarrow{\cal O}(U_0). 
 \end{equation*}
 Since ${\cal O}(U)=\bigcup_{i\geq0}{\cal O}(U_0)^{p^{-i}}$
  and ${\cal O}(V_0)$ is a finitely generated $k$-algebra, we may assume that $f^\ast({\cal O}(V_0))$ is contained in ${\cal O}(U_0)$ by replacing ${\cal O}(U_0)$ with ${\cal O}(U_0)^{p^{-i}}$ for a large $i$. Let $f_0^\ast\colon {\cal O}(V_0)\to {\cal O}(U_0)$ be the induced map. This is automatically a map of Hopf $k$-algebras, since ${\cal O}(U_0)\to {\cal O}(U)$ is injective and ${\cal O}(V_0)\to{\cal O}(U)$ is a map of Hopf $k$-algebras. Then we may take $f_0$ to be the map 
 $ U_0\to V_0$ which corresponds to $ {\cal O}(V_0)\to{\cal O}(U_0)$. 
   \qed}

We deduce some necessary results from this lemma. 

\begin{lm}\label{a}
Let $f\colon U\to V$ be a morphism of perfect unipotent commutative group schemes over $k$. Then the kernel of $f$ is represented by a perfect $k$-scheme. It 
has finitely many connected components and the connected component of the neutral element is perfect unipotent commutative. 
 \end{lm}
 \proof{
 Take  algebraic models $U_0,V_0,$ and $f_0\colon U_0\to V_0$ as in Lemma \ref{redalg}. Since the perfection commutes with fiber product, the perfection of $\mathop{\rm ker}(f_0)$ represents $\mathop{\rm ker}(f)$. The scheme $\mathop{\rm ker}(f)$ has finitely many connected components as $\mathop{\rm ker}(f_0)$ does. 
  
  Let $K^\circ$ be the neutral connected component of 
  $\mathop{\rm ker}(f_0)$. Let $K^\circ_{\rm red}$ be 
  the reduced closed subscheme of $K^\circ$ with the same underlying space as  $K^\circ$. This is smooth and unipotent. The kernel of ${\cal O}(K^\circ)\to {\cal O}(K^\circ_{\rm red})$ is a nilpotent ideal, which will be killed by the perfection. Therefore, the perfection of $K^\circ$ is isomorphic to that of $K^\circ_{\rm red}$. 
   The last assertion follows. 
  \qed}

\begin{lm}\label{morper}
Let $f\colon U\to V$ be a morphism of perfect unipotent commutative group $k$-schemes. Under the assumptions given below, we may find algebraic models $f_0\colon U_0\to V_0$ for $f$ with 
the following properties. 
\begin{enumerate}
\item Suppose that $f$ is injective on the underlying spaces. Then we may take $f_0$ to be a closed immersion.  Moreover, $f$ is a closed immersion. 
\item Suppose that $f$ is surjective on the underlying spaces and that $U,V$ have the same dimension. Then we may take $f_0$ to be finite \'etale. Moreover, $f$ is finite \'etale. 
\item Suppose that $f$ is bijective on the underlying spaces. Then we may take $f_0$ to be an isomorphism. 
Moreover, $f$ is an isomorphism. 
\end{enumerate}
\end{lm}
\proof{
Take algebraic models $U_0,V_0$ for $U,V,$ and $f_0\colon U_0\to V_0$ for $f$.  We modify our choices as follows. Let $K^\circ$ denote the connected component of the neutral element in ${\rm ker}(f_0)$. Note that $f_0$ is finite in either of the three cases we treat. This implies that  $K^\circ$ is a finite local group scheme, which will be killed after perfection. Thus we may replace $U_0$ with $U_0/K^\circ$ and $f_0$ with the induced map $U_0/K^\circ\to V_0$: we may further assume that $f_0$ is finite \'etale. 

$1.$ In this case, $\mathop{\rm ker}(f_0)$ is isomorphic to ${\rm Spec}(k)$, which implies that $f_0$ is a closed immersion. Since the perfection functor makes a closed immersion into a closed immersion, $f$ is a closed immersion, too. 

$2.$ As $f_0$ is \'etale, $U_0\times_{V_0}V$ is \'etale over a perfect scheme $V$. Therefore, $U_0\times_{V_0}V$ is perfect. We can show that this is isomorphic to the perfection of $U_0$. Therefore, $f$ is finite \'etale. 

$3.$ By $1.$ and $2.$, we know that $f$ is a closed immersion which is finite \'etale, which is equivalent to that $f$ is an isomorphism. Then we may obviously take $f_0$ to be an isomorphism. 
\qed}

 \begin{lm}\label{section}
 Let 
 \begin{equation*}
 0\to U\to V\to W\to 0
 \end{equation*}
 be a short exact sequence of perfect unipotent commutative group $k$-schemes. Then there exists a morphism $W\to V$ of $k$-schemes such that $W\to V\to W$ is the identity. 
 \end{lm}
\proof{
By Lemma \ref{morper}.1, we may find an algebraic model $f_0\colon U_0\to V_0$ for $f$ which is a closed immersion. Let $W_0:=V_0/U_0$, which is an affine smooth connected unipotent commutative group $k$-scheme. We view $V_0$ as a $U_0$-torsor on $W_0$ via the surjection $V_0\to W_0$.  
Note that any $U_0$-torsor on an affine scheme is trivial: this follows 
since any $\mathbb{G}_a$-torsor on an affine scheme is trivial and $U_0$ has a finite filtration given by closed group subschemes whose successive quotients are  $\mathbb{G}_a$. In particular, $V_0$ is trivial as a $U_0$-torsor on $W_0$, i.e., the projection $V_0\to W_0$ has a section. Consequently, the sequence 
\begin{equation*}
0\to U_0\to V_0\to W_0\to0
\end{equation*}
is exact 
 as a sequence of presheaves on $Sch/k$.  Then 
 the restriction of this sequence to $PSch/k$ remains exact. Therefore, we know that the perfection of $W_0$ is isomorphic to $W$. Then the projection $V\to W$ has a section as $V_0\to W_0$ does. 
 \qed}
\subsection{Duality for perfect unipotent commutative groups}
In this subsection, we review the duality theorem of perfect unipotent commutative group schemes, which is an  essential tool in this article. 
The dual group of a perfect unipotent group is commonly defined in terms of extensions by $\Qp/\Zp$ in the literature. 
In this article, we define it in terms of multiplicative sheaves since this is more useful for our purposes; it is known that the two definitions actually give the same dual groups, which is reviewed in Lemma \ref{compdual}. 

\subsection*{Reminders on multiplicative sheaves and extensions}
Let $\Lambda$ be a commutative ring as in Convention and Notation. Multiplicaitve $\Lambda$-sheaves are defined as follows. 
\begin{df}\label{mult}
Let $S$ be an $\F_p$-scheme and $G$ be a commutative group $S$-scheme. 
{\rm A multiplicative $\Lambda$-sheaf on $G/S$} is an invertible \'etale $\Lambda$-sheaf ${\cal L}$ on $G$ such that $m^\ast{\cal L}$ is isomorphic to ${\rm pr}_1^\ast{\cal L}\otimes_\Lambda
{\rm pr}_2^\ast{\cal L}$, where $m$ is the group law $G\times_SG\to G$ and ${\rm pr}_i$ denotes the $i$-th projection. 
\end{df}

We write ${\rm MS}_\Lambda(G/S)$ for the set of isomorphism classes of multiplicative $\Lambda$-sheaves on $G/S$, which we consider  as a group with the group law given by tensor product. 
For a morphism $T\to S$ of schemes, the pullback induces a group homomorphism ${\rm MS}_\Lambda(G/S)\to {\rm MS}_\Lambda(G_T/T)$, where $G_T=G\times_ST$: in this way, we usually view the groups of multiplicative sheaves as a presheaf on $Sch/S$ in what follows. 

We recall some basic results on multiplicative sheaves. 
\begin{lm}\label{basicmult}
Let $G$ be a commutative group $S$-scheme. 
\begin{enumerate}
\item Let $\cal L$ be a multiplicative $\Lambda$-sheaf on $G/S$. Let $0$ denote the $0$-section $S\to G$. Then the pullback $0^\ast{\cal L}$ is trivial. 
\item Let $\cal L$ be a multiplicative $\Lambda$-sheaf on $G/S$. Let $n\geq1$ be an integer which kills $G$. Then ${\cal L}^{\otimes n}$ is trivial. 
\item Suppose that $S$ is perfect. Let $G^{\rm perf}$ be the perfection of $G$, which is a group $S$-scheme. 
Then pulling back via $G^{\rm perf}\to G$ induces a bijection ${\rm MS}_\Lambda(G/S)\to{\rm MS}_\Lambda(G^{\rm perf}/S)$. 
\end{enumerate}
\end{lm}
\proof{Let $m$ denote the multiplication $G\times_SG\to G$. 

$1.$ Applying $0^\ast$ to $m^\ast{\cal L}\cong {\rm pr}_1^\ast{\cal L}\otimes_\Lambda
{\rm pr}_2^\ast{\cal L}$, the assertion follows. 

$2.$ Let $n_G$ denote the multiplication by $n$ on $G$. The multiplicativity of $\cal L$ implies that ${\cal L}^{\otimes n}$ is isomorphic to $n_G^\ast{\cal L}$. 
On the other hand, the map $n_G$ factors over the $0$-section $S\to G$. Then the assertion follows from $1.$ 

$3.$ The canonical map $G^{\rm perf}\to G$ induces an equivalence between the \'etale topoi. Similarly, 
the \'etale topoi of $G\times_SG$ and $G^{\rm perf}\times_SG^{\rm perf}$ are equivalent. Then the assertion follows. 
\qed}

Let $U$ be a commutative group $k$-scheme. 
We assume that $U$ is geometrically connected over $k$ and that it is killed by $p^n$ for some $n\geq1$. 
For a perfect $k$-scheme $S$, we let 
${\rm Ext}(U_S, \Qp/\Zp)$ be the group of extensions 
\[0\to\Qp/\Zp\to E\to U_S\to0, 
\]
where $\Qp/\Zp$ is viewed as a constant group scheme over $S$. In Lemma \ref{compdual} below, we recall the result that ${\rm MS}_\Lambda(U_S/S)$ can be identified with  ${\rm Ext}(U_S, \Qp/\Zp)$ once we fix an injective character $\Qp/\Zp\to\Lambda^\times$.

Let $S$ be a perfect $k$-scheme and let 
$E$ 
be an extension of $U_S$ by $\Qp/\Zp$ as above. For a  character $\psi\colon \Qp/\Zp\to\Lambda^\times$, the push forward of a $\Qp/\Zp$-torsor $E$ via
$\psi$ gives a $\Lambda^\times$-torsor on $U_S$. 
We write $E_\psi$ for the corresponding invertible \'etale $\Lambda$-sheaf on $U_S$. 
In \cite[Lemma 7.3]{Bo}, the following is proved at least when $S={\rm Spec}(k)$ and $U$ is a unipotent algebraic group. 
\begin{lm}\label{compdual}
Let the notation and the assumptions be as above. We suppose that $\psi$ is injective. Then the assignment $E\mapsto E_\psi$ constructed above induces a bijection 
\[{\rm Ext}(U_S,\Qp/\Zp)\cong {\rm MS}_\Lambda(U_S/S). 
\]
\end{lm}

\proof{
We view ${\rm MS}_\Lambda(U_S/S)$ as the subgroup of $H^1(U_S,\Lambda^\times)$ consisting of those invertible sheaves which have the multiplicative property. The assignment $E\mapsto E_\psi$ gives a group homomorphism ${\rm Ext}(U_S,\Qp/\Zp)\to 
H^1(U_S,\Lambda^\times)$.

Let $E$ be an extension of $U_S$ by $\Qp/\Zp$. We check that the element $E_\psi\in H^1(U_S,\Lambda^\times)$ is contained in ${\rm MS}_\Lambda(U_S/S)$. Indeed, from the commutativity of the following diagram 
\begin{equation*}
\xymatrix{
E\times_SE\ar[d]\ar[r]^-{m_E}&E\ar[d]\\
U_S\times_SU_S\ar[r]^-m&U_S, 
}
\end{equation*}
where $m_E$ is the multiplication of $E$, we have an equality ${\rm pr}_1^\ast E_\psi+{\rm pr}_2^\ast E_\psi=m^\ast E_\psi$ in $H^1(U_S,\Lambda^\times)$, which implies that $E_\psi$ is a multiplicative $\Lambda$-sheaf.  Thus $E\mapsto E_\psi$ defines a map ${\rm Ext}(U_S,\Qp/\Zp)\to {\rm MS}_\Lambda(U_S/S)$. 
We have to show that this is bijective. 
By Lemma \ref{ord}.1 below, it suffices to show that the composite map 
${\rm Ext}(U_S,p^{-n}\Zp/\Zp)\to{\rm Ext}(U_S,\Qp/\Zp)\to {\rm MS}_\Lambda(U_S/S)$ is bijective where $n\geq0$ is an integer such that $U$ is killed by $p^n$. 

We construct the inverse of this map. Let $\cal L$ be a multiplicative $\Lambda$-sheaf. By Lemma \ref{basicmult}.2., we know that ${\cal L}^{\otimes p^{n}}$ is trivial. The exact sequence 
\[0\to p^{-n}\Zp/\Zp\to\Lambda^\times\xrightarrow{p^n}\Lambda^\times\to0, 
\]
induces an exact sequence 
\[H^0(U_S,\Lambda^\times)\xrightarrow{p^n}
H^0(U_S,\Lambda^\times)\to H^1(U_S,p^{-n}\Zp/\Zp)\to H^1(U_S,\Lambda^\times)\xrightarrow{p^n}
H^1(U_S,\Lambda^\times). 
\]
Since $H^0(U_S,\Lambda^\times)\xrightarrow{p^n}
H^0(U_S,\Lambda^\times)$ is surjective, the map 
$H^1(U_S,p^{-n}\Zp/\Zp)\to H^1(U_S,\Lambda^\times)$ is injective. Consequently, there exists a unique $
p^{-n}\Zp/\Zp$-torsor $E$ on $U_S$ which maps to $\cal L$. Since $0^\ast{\cal L}$ is trivial by Lemma \ref{basicmult}.1., $0\colon S\to U_S$ lifts to a section $0_E\colon S\to E$. 
Moreover, the multiplicativity of $\cal L$ implies that the map $m\colon U_S\times_SU_S\to U_S$ extends to 
a map $m_E\colon E\times_SE\to E$ equivariantly to the map 
$p^{-n}\Zp/\Zp\times p^{-n}\Zp/\Zp\to p^{-n}\Zp/\Zp,(a,b)\mapsto a+b$. Replacing $m_E$ with $m_E+a$ for a suitable $a\in p^{-n}\Zp/\Zp$ if necessary, we may assume that $m_E(0_E,0_E)=0_E$. Then one can check that $m_E$ provides a structure of commutative group scheme to $E$ and that $E\to U_S$ is an extension by $\Qp/\Zp$ which maps to $\cal L$. 
\qed
}

The following is used in the proof of the lemma above. 
\begin{lm}\label{ord}
 Let $U$ be a commutative group $k$-scheme. Suppose that it is geometrically connected over $k$ and that there exists  
an integer $n\geq0$ such that $U$ is killed by $p^n$. 
 \begin{enumerate}
\item 
 For $S\in PSch/k$, the  map 
 \[{\rm Ext}_S(U_S,p^{-n}\Zp/\Zp)\to {\rm Ext}_S(U_S,\Qp/\Zp)
 \]
  induced from the inclusion $p^{-n}\Zp/\Zp\to \Qp/\Zp$ is an isomorphism. 
  \item Fix an injective character $\psi\colon\Qp/\Zp\to\Lambda^\times$. For $S\in PSch/k$, the composite map ${\rm Ext}_S(U_S,p^{-n}\Zp/\Zp)\to {\rm Ext}_S(U_S,\Qp/\Zp)\xrightarrow{\text{Lemma \ref{compdual}}}{\rm MS}_\Lambda(U_S/S)$ only depends on the restriction $\psi|_{p^{-n}\Zp/\Zp}$. 

  \end{enumerate}
 \end{lm}
 \proof{
 $1.$ 
 The exact sequence 
 \[0\to p^{-n}\Zp/\Zp\to\Qp/\Zp\xrightarrow{p^n}\Qp/\Zp\to0
 \]
  induces a long exact  sequence
 \begin{align*}
 {\rm Hom}_S(U_S,\Qp/\Zp)\to{\rm Ext}_S(U_S,p^{-n}\Zp/\Zp)\to{\rm Ext}_S(U_S,\Qp/\Zp)\xrightarrow{p^n}{\rm Ext}_S(U_S,\Qp/\Zp). 
 \end{align*}
 The group ${\rm Hom}_S(U_S,\Qp/\Zp)$ is trivial since $U$ is geometrically connected over $k$ while $\Qp/\Zp$ is constant.  Moreover, the map 
 ${\rm Ext}_S(U_S,\Qp/\Zp)\xrightarrow{p^n}{\rm Ext}_S(U_S,\Qp/\Zp)$ is zero by the assumption $p^n U=0$. Therefore, we have ${\rm Ext}_S(U_S,p^{-n}\Zp/\Zp)\cong{\rm Ext}_S(U_S,\Qp/\Zp)$, as claimed. 
 
 $2.$ The assertion follows as the composite map is given by pushing forward via $\psi|_{p^{-n}\Zp/\Zp}$, which is defined similarly as the isomorphism in  Lemma \ref{compdual}. 
 \qed}

 We recall the following vanishing result. 
 \begin{lm}\label{vancoh}
 Suppose that $k$ is algebraically closed. 
 Let $U$ be a perfect unipotent commutative group $k$-scheme. Let $\cal L$ be a multiplicative $\Lambda$-sheaf on $U$ which is non-trivial. Then $R\Gamma_c(U,{\cal L})$ is acyclic. 
 \end{lm}
 \proof{
 Take an algebraic model $U_0$ for $U$. We may identify $\cal L$ with an \'etale sheaf on $U_0$ and $R\Gamma_c(U,{\cal L})$ with $R\Gamma_c(U_0,{\cal L})$. For algebraic groups, the same statement is proved in  \cite[Lemma 9.4]{Bo}. 
 \qed
 }
 
\subsection*{Duality}

We define the dual group as follows. 
\begin{df}
Let $U$ be a perfect unipotent commutative group scheme over $k$. We define the dual group $U^\ast_\Lambda$ of $U$ to be the  presheaf on $PSch/k$ given by 
$U^\ast_\Lambda(S):={\rm MS}_\Lambda(U_S/S)$, where $U_S$ is the base change $U\times_kS$. 
\end{df}
We simply write $U^\ast$ for $U^\ast_\Lambda$ if no confusion occurs. By Lemma \ref{compdual}, the presheaf $U^\ast$ can be identified with $S\mapsto {\rm Ext}(U_S,\Qp/\Zp)$ when one fixes an injective character $\Qp/\Zp\to\Lambda^\times$. 
 From the following theorem, we know that the assignment $U\mapsto U^\ast$ gives a contravariant functor from the category of perfect unipotent commutative group $k$-schemes to itself.

\begin{thm}\label{dualp}
\begin{enumerate}
\item{(\cite[Proposition 1.2.1]{Beg})} The dual group $U^\ast$ of a perfect unipotent commutative group $k$-scheme $U$  is representable by a perfect unipotent commutative group $k$-scheme. Moreover, we have $ \dim U^\ast=\mathop{\rm dim}U$. 
\item{(\cite[Proposition 1.2.1]{Beg})} Let 
\begin{equation*}
0\to U\to V\to W\to0
\end{equation*}
be a short exact sequence of perfect unipotent commutative group $k$-schemes. Then the induced sequence 
\begin{equation*}
0\to W^\ast\to V^\ast\to U^\ast\to0
\end{equation*}
is again exact. 
\item Let $f\colon U\to V$ be a morphism of perfect unipotent commutative group $k$-schemes which is finite \'etale. Then $f^\ast\colon V^\ast\to U^\ast$ is finite \'etale. Moreover, we have $\deg f=\deg f^\ast$. 
\end{enumerate}
\end{thm}
\proof{ By fixing an injective character $\Qp/\Zp\to\Lambda^\times$, we identify $U^\ast(S)$ and ${\rm Ext}_S(U_S,\Qp/\Zp)$. In \cite[Proposition 1.2.1]{Beg}, it is proved that $U^\ast$ is representable by a perfect unipotent commutative group scheme and that the duality functor is exact\footnote{In loc.~cit., a perfect unipotent commutative group is called as groupe quasi-alg\'ebrique unipotent connexe. }. Moreover, $U^\ast$ is isogeneous to $U$. In particular, we have $\dim U^\ast=\dim U$.  Hence the assertions $1.$~and $2.$~follow.

$3.$ Let $K$ be the kernel of $f$, which is a finite \'etale group scheme over $k$. Applying $R{\rm Hom}(-,\Qp/\Zp)$ to $0\to K_S\to U_S\to V_S\to0$, 
we get a long exact sequence 
\begin{align*}
{\rm Hom}(U_S,\Qp/\Zp)\to{\rm Hom}&(K_S,\Qp/\Zp)
\\&\to {\rm Ext}(V_S,\Qp/\Zp)\xrightarrow{f^\ast}{\rm Ext}(U_S,\Qp/\Zp)\to{\rm Ext}(K_S,\Qp/\Zp). 
\end{align*}
 As 
$\Qp/\Zp$ is an injective abelian group, the sheafification of $S\mapsto {\rm Ext}(K_S,\Qp/\Zp)$ vanishes. Thus $f^\ast$ is surjective. Since we have $\dim V^\ast=\dim V=\dim U=\dim U^\ast$, 
 Lemma \ref{morper}.2 implies that $f^\ast$ is finite \'etale. 
As $U/k$ is geometrically connected, we have ${\rm Hom}(U_S,\Qp/\Zp)=0$, which implies the assertion on the degree. 
\qed}

\subsection{Multiplicative sheaves on $\G_{a}^{\rm perf}$}\label{multGa}
 Let $\mathbb{G}_{a,k}^{\rm perf}$ be the perfection of $\mathbb{G}_{a,k}$. In this subsection, we  review  the concrete description of $\mathbb{G}_{a,k}^{{\rm perf}\ast}$ in terms of Artin--Schreier theory.

 Let $\psi\colon \F_p\to \Lambda^\times$ be a non-trivial additive character and let ${\cal L}_\psi$ be the Artin--Schreier sheaf on $\G_{a,k}$ defined from $\psi$ 
 and the Artin--Scheier covering $\G_{a,k}\to \G_{a,k}, x\mapsto x^p-x$. 
 Let $S$ be a perfect $k$-scheme. For $f\in\Gamma(S,{\cal O}_S)$, we let ${\cal L}_\psi(fx)$ be the pullback of ${\cal L}_\psi$ via the map $\G_{a,S}\to\G_{a,S}, x\mapsto fx$. We also write ${\cal L}_\psi(fx)$ for its pullback to $\G_{a,S}^{\rm perf}$. 
 
 \begin{pr}\label{LGa}(cf. \cite[Proposition 1.20]{Sa}) 
 Let $S$ be a perfect $k$-scheme. Then the assignment $f\mapsto {\cal L}_\psi(fx)$ gives an isomorphism $\Gamma(S,{\cal O}_S)\to \G_{a,k}^{{\rm perf}\ast}(S)$. 
 \end{pr}
 \proof{
 Let 
 \begin{equation*}
 0\to\F_p\to\G_{a,S}\xrightarrow{x\mapsto x^p-x}\G_{a,S}\to 0
 \end{equation*}
 be the Artin--Schreier extension and $[AS]$ denote the element of ${\rm Ext}(\G_{a,S},\F_p)$ corresponding to this extension. For an element $f\in \Gamma(S,{\cal O}_S)$, let $f^\ast[AS]$ be the extension class obtained by the pullback of $[AS]$ via $\G_{a,S}\to\G_{a,S},x\mapsto fx$. In \cite[Proposition 1.20]{Sa}, it is proved that the assignment $f\mapsto f^\ast[AS]$ induces a bijection $\Gamma(S,{\cal O}_S)\to {\rm Ext}(
 \G_{a,S},\F_p)$: indeed, it is proved in loc.~cit.~with the extra assumption that $S$ is quasi-compact. The general case follows from the quasi-compact case as the both of the source and the target 
 define Zariski sheaves on $S$.

 By Lemma \ref{compdual} and Lemma \ref{ord}.2, the character $\psi$ induces an isomorphism 
 \[{\rm Ext}_S(\G_{a,S},\F_p)\cong {\rm MS}_\Lambda(\G_{a,S}/S).\] 
 Here we identify $\F_p$ and $p^{-1}\Zp/\Zp$ via $1\mapsto 1/p.$ Moreover, by Lemma \ref{basicmult}.3., ${\rm MS}_\Lambda(\G_{a,S}/S)$ is canonically isomorphic to $\G_{a,k}^{{\rm perf}\ast}(S)$. 
 Then it remains to  check that the composite map $\Gamma(S,{\cal O}_S)\to \G_{a,k}^{{\rm perf}\ast}(S)$ is given by $f\mapsto {\cal L}_\psi(fx)$, which is straightforward. 
 \qed}

 \subsection{Bimultiplicative sheaves}

 \begin{df}\label{bimult}
 Let $S$ be an $\F_p$-scheme and $G,H$ be commutative group schemes over $S$. By {\rm a bimultiplicative $\Lambda$-sheaf on $G\times_SH$}, we mean an \'etale invertible $\Lambda$-sheaf $\cal B$ on $G\times_SH$ with the following property:
 when $G\times_SH$ is viewed as a group scheme over $G$ (resp. over $H$) 
 via the first (resp. second) projection, it is a multiplicative sheaf in the sense of Definition \ref{mult}. 
 \end{df}

Let $U,V$ be perfect unipotent commutative group schemes over  $k$. Then a bimultiplicative sheaf $\cal B$ on $U\times_kV$ induces a morphism of group $k$-schemes $l_{\cal B}\colon U\to V^\ast$ defined by  $f\mapsto (f\times {\rm id}_V)^\ast{\cal B}$. We often write  ${\cal B}(f,-)$ for $(f\times{\rm id}_V )^\ast{\cal B}$ to ease the notation. 
The assignment ${\cal B}\mapsto l_{\cal B}$ induces a bijection between the set of bimultiplicative $\Lambda$-sheaves on $U\times_kV$ and that of group homomorphisms $U\to V^\ast$. More generally, for a perfect $k$-scheme $S$, 
the same construction gives an identification between the set of bimultiplicative $\Lambda$-sheaves on $U_S\times_SV_S$ and that of group homomorphisms $U_S\to V^\ast_S$.

Similarly, one can define $r_{\cal B}\colon V\to U^\ast$ by setting $r_{\cal B}(f):=({\rm id}_U\times f)^\ast{\cal B}.$

\begin{df}\label{dymbimdf}
Let $U$ be a commutative group scheme over an $\F_p$-scheme $S$. Let $\cal B$ be a bimultiplicative $\Lambda$-sheaf on $U\times_SU$. 
\begin{enumerate}
\item We say that $\cal B$ is {\rm symmetric} if $\tau^\ast{\cal B}$ is isomorphic to $\cal B$, where $\tau\colon U\times_S U\to U\times_SU$ is the map $(x,y)\mapsto (y,x)$. 
\item Further suppose that $S={\rm Spec}(k)$ for a perfect field $k$ and that $U$ is perfect unipotent. We say that $\cal B$ is {\rm isogeneous} if it is symmetric and if the map $l_{\cal B}\colon U\to U^\ast$ is finite \'etale. If ${\cal B}$ is isogeneous, then we define the degree of $\cal B$ to be the degree of $l_{\cal B}$. 
Note that the degree is always a power of $p$. 
\item With the same assumption as $2.$, we say that $\cal B$ is {\rm non-degenerate} if it is isogeneous of degree $1.$ 
\end{enumerate}
\end{df}
When $S={\rm Spec}(k)$ and $U$ is perfect unipotent, the condition that $\cal B$ is symmetric is equivalent to requiring that $r_{\cal B}=l_{\cal B}$.

\subsection{Bimultiplicative sheaves on $\G_a^{\rm perf}$}

In this subsection, we give a classification result for symmetric bimultiplicative $\Lambda$-sheaves on $\G_a^{\rm perf}\times\G_a^{\rm perf}$. 

Fix a non-trivial character $\psi\colon\F_p\to\Lambda^\times$ and let  
${\cal L}_\psi$ be the Artin--Schreier sheaf on $\mathbb{G}_{a,\F_p}$ defined in the same way as in \ref{multGa}. 
We identify $\G_{a}^{{\rm perf}\ast}$ with $\G_a^{\rm perf}$ via $f\mapsto {\cal L}_\psi(fx)$ (Proposition \ref{LGa}). 
Let $A$ be a perfect $\F_p$-algebra. By an additive function on $\G_{a,A}^{\rm perf}$, we mean an element in $A[X,X^{p^{-1}},X^{p^{-2}},\dots]$ which is of the form $\sum_{i\in\mathbb{Z}}a_iX^{p^{i}}$. The additive functions are naturally identified with the endomorphisms $\G_{a,A}^{\rm perf}\to \G_{a,A}^{\rm perf}$ as a group $A$-scheme.

\begin{lm}\label{clB}
Let $A$ be a perfect $\F_p$-algebra and let $\cal B$ be a symmetric bimultiplicative $\Lambda$-sheaf on $\G_{a,A}^{\rm perf}\times_A\G_{a,A}^{\rm perf}$. Then there exists an additive function $f(X)$ which is of the form $a_0X+\sum_{i\geq1}(a_iX^{p^i}+a_i^{p^{-i}}X^{p^{-i}})$ 
such that ${\cal B}$ is isomorphic to ${\cal L}_\psi(f(x)y)$ where $x,y$ are the coordinates of $\G_{a,A}^{\rm perf}\times_A\G_{a,A}^{\rm perf}$.
\end{lm}
\proof{
By the assignment ${\cal B}\mapsto r_{\cal B}$, we identify the bimultiplicative $\Lambda$-sheaves on $\G_{a,A}^{\rm perf}\times_A\G_{a,A}^{\rm perf}$ with the endomorphisms of $\G_{a,A}^{\rm perf}$, hence with the additive functions: namely, for a bimultiplicative $\Lambda$-sheaf $\cal B$, there exists a unique $f(X)=\sum_{i\in\mathbb{Z}}a_iX^{p^{i}}$ such that ${\cal B}\cong {\cal L}_\psi(f(x)y)$. Therefore, what we have to check is that $\cal B$ is symmetric if and only if $a_i^{p^{-i}}=a_{-i}$ for $i\geq1$. 

First suppose that $f(X)$ satisfies $a_i^{p^{-i}}=a_{-i}$ for $i\geq1$. Then we have 
\begin{equation*}
f(X)Y-Xf(Y)=\sum_{i\geq1}((a_iX^{p^i}+a_i^{p^{-i}}X^{p^{-i}})Y-X(a_iY^{p^i}+a_i^{p^{-i}}Y^{p^{-i}})). 
\end{equation*}
It suffices to show that, for each $i\geq1$, there exists $g\in A[X,Y,X^{p^{-1}},Y^{p^{-1}},\dots]$ satisfying $(a_iX^{p^i}+a_i^{p^{-i}}X^{p^{-i}})Y-X(a_iY^{p^i}+a_i^{p^{-i}}Y^{p^{-i}})=g^p-g$. By a simple  computation, we know that the \'etale sheaf ${\cal L}_\psi((x^{p^i}-a_ix)(y^{p^i}-a_iy))$ on $\G_{a,A}^{\rm perf}\times_A\G_{a,A}^{\rm perf}$ is isomorphic to ${\cal L}_\psi(-(a_ix^{p^i}-(1+a_i^2)x+a_i^{p^{-i}}x^{p^{-i}})y)$. Since this is clearly symmetric, 
\begin{align*}
-(a_iX^{p^i}-(1+a_i^2)X+a_i^{p^{-i}}X^{p^{-i}})Y+&(a_iY^{p^i}-(1+a_i^2)Y+a_i^{p^{-i}}Y^{p^{-i}})X\\&=-(a_iX^{p^i}+a_i^{p^{-i}}X^{p^{-i}})Y+(a_iY^{p^i}+a_i^{p^{-i}}Y^{p^{-i}})X
\end{align*}
 must be of the form $h^p-h$. Then we can take $g=-h$. 
Explicitly, $g$ can be given by $g=\sum_{j=0}^{i-1}(a_i^{p^{-i}}XY^{p^{-i}}-a_i^{p^{-i}}X^{p^{-i}}Y)^{p^j}$. 

Next suppose that ${\cal B}$ is symmetric. This means that there exists an element $g\in A[X,Y,X^{p^{-1}},Y^{p^{-1}},\dots]$ such that 
\begin{equation}\label{A}
f(X)Y-Xf(Y)=g^p-g. 
\end{equation}
From the computation above, there exists $h_{i,a}\in A[X,Y,X^{p^{-1}},Y^{p^{-1}},\dots]$ for each $i\geq1$ and $a\in A$ such that $(aX^{p^i}+a^{p^{-i}}X^{p^{-i}})Y-X(aY^{p^i}+a^{p^{-i}}Y^{p^{-i}})=h_{i,a}^p-h_{i,a}$. Then replacing $f(X)$ with $f(X)-a_0X-\sum_{i\geq1}(a_{-i}^{p^i}X^{p^i}+a_{-i}X^{p^{-i}})$, we may assume that $a_i=0$ for $i\leq0$. In particular, $f(X)$ is contained in $A[X]$. In this case, we may pick up $g$ in (\ref{A}) from the polynomial ring $A[X,Y]$ (this is a consequence of the invariance of \'etale cohomology under the Frobenius pullback). Let us write $g=\sum_{n\geq0}g_n(X)Y^n$ where $g_n(X)\in A[X]$. Then (\ref{A}) can be rewritten as 
\begin{equation}\label{x}
f(X)Y-\sum_{i\geq1}a_iXY^{p^i}=\sum_{n\geq0}g_n(X)^pY^{np}-\sum_{n\geq0}g_n(X)Y^n. 
\end{equation}
Suppose that $a_i\neq0$ for some $i\geq1$. Let $m$ be the largest integer $\geq1$ such that the coefficient 
of $Y^m$ in (\ref{x}) is non-zero. From the left-hand side, we know that the coefficient of $Y^m$ is of the form $aX$ for some $a\in A\setminus\{0\}$ while we know from the right-hand side that it is of the form $h(X)^p$ for some $h(X)\in A[X]$, a contradiction. Thus we must have $f(X)=0$, which shows the assertion. 
\qed}

\subsection{Some descent results for invertible \'etale $\Lambda$-sheaves}\label{dessubs}
In this preliminary subsection, we collect some  descent results for invertible \'etale sheaves which are not equipped with rigidification. 

We give notation. For a topological space $T$, we write $LC(T,\Lambda^\times)$ for the group of locally constant functions $T\to \Lambda^\times$. We repeatedly use the following fact: let $X$ be an scheme and $\cal F$ be an invertible \'etale $\Lambda$-sheaf on it. Then the automorphism group ${\rm Aut}({\cal F})$ is canonically identified with $LC(|X|,\Lambda^\times)$, where $|X|$ is the underlying space.

 \begin{lm}\label{des}
 Let $f\colon X\to Y$ be a morphism of schemes which is universally open, surjective, and geometrically connected. Let $\cal F$ be an invertible \'etale $\Lambda$-sheaf on $X$. Suppose that there exists an isomorphism ${\rm pr}_1^\ast{\cal F}\cong {\rm pr}_2^\ast{\cal F}$ where ${\rm pr}_i\colon X\times_YX\to X$ denotes the $i$-th projection. Then there exists an isomorphism $\sigma\colon {\rm pr}_1^\ast{\cal F}\cong {\rm pr}_2^\ast{\cal F}$ which satisfies the cocycle condition. Consequently, there exists an invertible $\Lambda$-sheaf $\cal G$ on $Y$ such that $f^\ast{\cal G}\cong\cal F$. 
 \end{lm}
 \proof{
 By the assumption on $X\to Y$ , the composition with  $f$ gives a bijection between $LC(|Y|,\Lambda^\times)$ and $LC(|X|,\Lambda^\times)$. Moreover, 
 the sets $
 LC(|Y|,\Lambda^\times)$, $LC(|X|,\Lambda^\times)$, and $ LC(|X\times_YX\times_YX|,\Lambda^\times)$ are  canonically identified: the important point is that the  maps $
 {\rm pr}_i^\ast\colon LC(|X|,\Lambda^\times)\to LC(|X\times_YX\times_YX|,\Lambda^\times)$ are equal to each other as they are the same as the composition of 
 \begin{equation*}
 LC(|X|,\Lambda^\times)\xrightarrow{(f^\ast)^{-1}}LC(|Y|,\Lambda^\times)\to LC(|X\times_YX\times_YX|,\Lambda^\times). 
 \end{equation*}
 
 Choose any isomorphism $\sigma'\colon{\rm pr}_1^\ast{\cal F}\xrightarrow{\cong} {\rm pr}_2^\ast{\cal F}$ and consider the composite map 
 \begin{equation*}
 {\rm pr}_{1,2}^\ast\sigma'\circ{\rm pr}_{2,3}^\ast\sigma'\circ({\rm pr}_{1,3}^\ast\sigma')^{-1}\colon {\rm pr}_1^\ast{\cal F}\to {\rm pr}_1^\ast{\cal F}
 \end{equation*}
  of sheaves on $X\times_YX\times_YX$. As ${\rm pr}_1^\ast{\cal F}$ is locally constant of rank $1$, 
  this map is given by the multiplication by some 
  $\alpha\in LC(|X\times_YX\times_YX|,\Lambda^\times)$. Let $\beta\colon|X|\to\Lambda^\times$ be a locally constant map such that $|X\times_YX\times_YX|\xrightarrow{{\rm pr}_i}|X|\xrightarrow{\beta}\Lambda^\times$ is equal to $\alpha$. Then $\beta^{-1}\sigma'$ satisfies the cocycle condition. 
  \qed}

Let $S$ be a scheme and $X$ be an $S$-scheme with a specified section $\sigma\colon S\to X$. We write $LS^{\rk=1}(X/S,\sigma)$ for the group of 
isomorphism classes of those invertible \'etale $\Lambda$-sheaves $\cal F$ whose pullback $\sigma^\ast{\cal F}$ is isomorphic to $\Lambda$. 
On the other hand, we write $LS^{\rk=1}_1(X/S,\sigma)$ for the group of isomorphism classes of pairs $({\cal F},\alpha)$ where 
$\cal F$ is an invertible \'etale $\Lambda$-sheaf on $X$ and $\alpha$ is an isomorphism $\sigma^\ast{\cal F}\cong \Lambda$. We have an obvious map $LS^{\rk=1}_1(X/S,\sigma)\to  LS^{\rk=1}(X/S,\sigma)$ by forgetting $\alpha$. 
\begin{lm}\label{rigid}
We further assume that $X\to S$ is universally open and geometrically connected. 
\begin{enumerate}
\item For a pair $({\cal F},\alpha)$ of an invertible \'etale $\Lambda$-sheaf on $X$ and an isomorphism 
$\alpha\colon\sigma^\ast{\cal F}\cong \Lambda$, its  automorphism group is trivial. 
\item The presheaf on $Sch/S$ given by $S'\mapsto
LS^{\rk=1}_1(X_{S'}/S',\sigma_{S'}), $ where $X_{S'}, \sigma_{S'}$ are the base changes, is an \'etale sheaf. 
\item The map $LS^{\rk=1}_1(X/S,\sigma)\to  LS^{\rk=1}(X/S,\sigma)$ is an isomorphism. 
\item The presheaf on $Sch/S$ given by $S'\mapsto
LS^{\rk=1}(X_{S'}/S',\sigma_{S'})$ is an \'etale sheaf. 
\end{enumerate}
\end{lm}
\proof{
Let $f$ denote the structure map $X\to S$. 
By the assumption on $f$, the map $LC(|S|,\Lambda^\times)\to LC(|X|,\Lambda^\times)$ given by $g\mapsto g\circ f$ is bijective. Consequently, the map 
$LC(|X|,\Lambda^\times)\to LC(|S|,\Lambda^\times)$ given by $g\mapsto g\circ\sigma$ is also bijective 
as the composite map $LC(|S|,\Lambda^\times)\to LC(|X|,\Lambda^\times)\to LC(|S|,\Lambda^\times)$ is equal to the identity.

$1.$ Let $g\colon{ \cal F}\to{\cal F}$ be an automorphism which is compatible with the rigidification $\alpha$. The isomorphism $g$ is given by the multiplication by a certain element in $LC(|X|,\Lambda^\times)$, which we also write $g$. The condition that $g$ is compatible with $\alpha$ is equivalent to that $g\circ\sigma=1$. Therefore, we have $g=1$ as desired. 

$2.$ This follows from the assertion $1.$ 

$3.$ The map is clearly surjective. We show the injectivity. Let $\cal F$ be an invertible \'etale $\Lambda$-sheaf on $X$ and $\alpha,\alpha'$ be isomorphisms $\sigma^\ast{\cal F}\cong\Lambda$. 
For such data, we have to find a locally constant function $g\colon |X|\to\Lambda^\times$ which makes the diagram
\begin{equation*}
\xymatrix{
\sigma^\ast{\cal F}\ar[rr]^-{g\circ\sigma}\ar[rd]_-\alpha
&&\sigma^\ast{\cal F}\ar[ld]^-{\alpha'}
\\
&\Lambda&
}
\end{equation*}
commutative. Such a $g$ does exist as the map $LC(|X|,\Lambda^\times)\to LC(|S|,\Lambda^\times), h\mapsto h\circ\sigma$ is bijective. 

$4.$ This is a consequence of $2.$~and $3.$
\qed}

\begin{cor}\label{fedes}
 Let $S$ be an $\F_p$-scheme and let $G$ be an     $S$-scheme which is universally open and geometrically connected over $S$. 
 Let $\cal M$ be an invertible \'etale $\Lambda$-sheaf on $G$. 
 Let $a\colon S'\to S$ be a  morphism of schemes which is smooth surjective. 
  For an object $F$ over $S$, we write $F_{S'}$ for the base change to $S'$. 

\begin{enumerate}
\item Let  $\sigma$ be a section $S\to G$. 
Suppose that the locally constant sheaf ${\cal M}_{S'}$ on $G_{S'}$ is trivial and that $\sigma^\ast\cal M$ is trivial. Then $\cal M$ is trivial. 
\item
Assume that $G$ is equipped with a structure of a commutative group $S$-scheme with the $0$-section denoted by $0$. 
If ${\cal M}_{S'}$ is a multiplicative $\Lambda$-sheaf on $G_{S'}$ and $0^\ast{\cal M}$ is trivial, then $\cal M$ is multiplicative. 
\end{enumerate}
 \end{cor}
 \proof{
 A smooth surjective map has a section \'etale locally. Thus 
 we may assume that $a$ is \'etale. 
 
 $1.$ This is a consequence of Lemma \ref{rigid}.4. 
 
 $2.$ 
 Set ${\cal M}':=m^\ast{\cal M}\otimes{\rm pr}_1^\ast{\cal M}^{-1}\otimes{\rm pr}_2^\ast{\cal M}^{-1}$, where $m$ denotes the map $G\times_SG\to G, (x,y)\mapsto x+y$. Its pullback via the section $\sigma:=(0,0)\colon S\to G\times_SG$ is trivial. Then Lemma \ref{rigid}.4.~implies that  ${\cal M}'$ is trivial. 
 \qed}

\section{Quadratic sheaves: definitions and basic properties}

Let $k$ be a perfect field of characteristic $p>0$  and $\Lambda$ be 
 our coefficient ring which has the property  indicated in Convention and Notation. Quadratic sheaves are defined as follows. 
\begin{df}\label{quaddef}
Let $S$ be an $\F_p$-scheme and $U$ be a commutative group $S$-scheme. 
\begin{enumerate}
\item 
{\rm A quadratic $\Lambda$-sheaf} on $U/S$ is an invertible \'etale $\Lambda$-sheaf $\cal Q$ on $U$ with the following property: set ${\cal B}_{\cal Q}:=m^\ast{\cal Q}\otimes{\rm pr}_1^\ast{\cal Q}^{-1}\otimes{\rm pr}_2^\ast{\cal Q}^{-1}$. This is a bimultiplicative $\Lambda$-sheaf on $U\times_SU$ in the sense of Definition \ref{bimult}. 
\item Suppose that $S={\rm Spec}(k)$ and that $U$ is a perfect unipotent commutative group $k$-scheme. 
A quadratic $\Lambda$-sheaf $\cal Q$ on $U$ is said  to be  {\rm isogeneous of degree $p^n$ } if the associated bimultiplicative $\Lambda$-sheaf ${\cal B}_{\cal Q}$ is isogeneous of degree $p^n$ in the sense of Definition \ref{dymbimdf}.2. We write $\deg{\cal Q}$ for the degree. 
\item With the same assumption as $2.$, we say that 
$\cal Q$ is {\rm non-degenerate} if it is isogeneous of degree $1.$ 
\end{enumerate}
\end{df}
For a quadratic $\Lambda$-sheaf $\cal Q$, we usually write $l_{\cal Q}$ for $l_{{\cal B}_{\cal Q}}$ for simplicity.

We prove some basic properties of quadratic sheaves. 
\begin{lm}\label{pullQ}
Let $S$ be an $\F_p$-scheme and $f\colon U\to V$ be a morphism of commutative group $S$-schemes. Let $\cal Q$ be a quadratic $\Lambda$-sheaf on $V$. 
\begin{enumerate}
\item The pullback $f^\ast\cal Q$ is quadratic. 
\item Further suppose that $S={\rm Spec}(k)$ and that $f$ is a finite \'etale morphism between perfect unipotent commutative group $k$-schemes. We also suppose that $\cal Q$ is isogeneous. Then $f^\ast{\cal Q}$ is isogeneous and its degree is equal to $\deg{\cal Q}\cdot(\deg f)^2$. 
\end{enumerate}
\end{lm}
\proof{
$1.$ 
It follows since ${\cal B}_{f^\ast{\cal Q}}$ is isomorphic to $(f\times f)^\ast{\cal B}_{\cal Q}$. 

$2.$ The map $l_{f^\ast{\cal Q}}$ is equal to the composite map 
\[
U\to V\xrightarrow{l_{\cal Q}}V^\ast\xrightarrow{f^\ast}U^\ast.
\]
The assertion follows as we have $\deg f=\deg f^\ast$ (Theorem \ref{dualp}.3). 
\qed}

\begin{lm}\label{order}
Let $\cal Q$ be a quadratic $\Lambda$-sheaf on a perfect unipotent commutative group $k$-scheme $U$. 
\begin{enumerate}
\item Let $0\in U(k)$ be the unit element. Then $0^\ast{\cal Q}$ is trivial. 
\item 
Let $n\geq1$ be an integer such that $p^n$ kills $U$. Then ${\cal Q}^{\otimes p^{2n}}$ is trivial. If $p$ is odd, then ${\cal Q}^{\otimes p^n}$ is trivial. 
\end{enumerate}
\end{lm}
\proof{
For a morphism $f\colon T\to U$, we write ${\cal Q}(f)$ for $f^\ast{\cal Q}$. 

$1.$ It follows since ${\cal Q}(0)^{-1}\cong {\cal B}_{\cal Q}(0,0)$ is trivial by Lemma \ref{basicmult}.1. 

$2.$ We show that ${\cal Q}^{\otimes p^{2n}}$ is trivial. Inductively on $j$, we find an isomorphism of invertible sheaves on $U$
\begin{equation*}
{\cal Q}(j\cdot{\rm id}_U)\otimes{\cal Q}^{\otimes(-j)}\cong \bigotimes_{i=1}^{j-1}{\cal B}_{\cal Q}(i\cdot{\rm id}_U,{\rm id}_U).
\end{equation*}
Consequently, we know that ${\cal Q}^{\otimes(-p^n)}\cong {\cal Q}(p^n\cdot{\rm id}_U)\otimes{\cal Q}^{\otimes(-p^n)}$ is a product of sections of ${\cal B}_{\cal Q}$. 
Then the assertion follows from the fact that ${\cal B}_{\cal Q}$ is killed by $p^n$, which is a consequence of 
Lemma \ref{basicmult}.2. 

We show that ${\cal Q}^{\otimes p^n}$ is actually trivial when $p$ is odd. Let ${\cal Q}':={\cal B}_{{\cal Q}}(\frac{1}{2}\cdot{\rm id}_U, {\rm id}_U)$, which is a quadratic sheaf with ${\cal B}_{{\cal Q}'}={\cal B}_{{\cal Q}}$. Thus ${\cal Q}\otimes{\cal Q}'^{-1}$ is a multiplicative $\Lambda$-sheaf on $U$. Then the assertion follows as ${\cal B}_{\cal Q}$ and ${\cal Q}\otimes{\cal Q}'^{-1}$ are killed by $p^n$ (Lemma \ref{basicmult}.2). 
\qed}

The following results are consequences of Corollary \ref{fedes}.

 \begin{lm}\label{desmq}
 Let $S$ be an $\F_p$-scheme and let $U,V$ be  commutative group $S$-schemes which are universally open and geometrically connected. Let 
$f\colon U\to V$ be a morphism of $S$-group schemes. 
 Let $\cal L$ be an invertible \'etale $\Lambda$-sheaf on $V$. Suppose that $f$ is smooth surjective. Then the following hold true. 
 \begin{enumerate}
 \item Suppose that $f^\ast\cal L$ is multiplicative. Then $\cal L$ is multiplicative. 
 \item Suppose that $f^\ast\cal L$ is quadratic. Then $\cal L$ is quadratic. 
 \end{enumerate}
 \end{lm}
 \proof{
 Let $0\colon S\to V$ be the $0$-section. Note that, in the both cases, the restriction $0^\ast{\cal L}$ is isomorphic to $\Lambda$ as $0$ factors through the $0$-section of $U$.

 $1.$ Let ${\cal M}:=m^\ast{\cal L}\otimes{\rm pr}_1^\ast{\cal L}^{-1}\otimes
{\rm pr}_2^\ast{\cal L}^{-1}$, where $m\colon V\times_S V\to V,(x,y)\mapsto x+y$.  We have to show that this is trivial, provided that $(f\times f)^\ast{\cal M}$ is trivial. 
Consider the following sequence of maps 
\begin{equation*}
U\times_SU\xrightarrow{f\times{\rm id}}V\times_SU\xrightarrow{{\rm id}\times f}V\times_S V. 
\end{equation*}
Applying Corollary \ref{fedes}.1 to $(S,G,\sigma)=(V,V\times_SU,V\to V\times_SU;x\mapsto(x,0))$, we know that $({\rm id}\times f)^\ast{\cal M}$ is trivial. Then, applying Corollary \ref{fedes}.1 again to $(S,X,\sigma)=(V,V\times_SV\xrightarrow{{\rm pr}_2}V,V\to V\times_SV;x\mapsto(0,x))$, we know that $\cal M$ is trivial, hence the assertion.

 $2.$ Let ${\cal M}:=m^\ast{\cal L}\otimes{\rm pr}_1^\ast{\cal L}^{-1}\otimes
{\rm pr}_2^\ast{\cal L}^{-1}$. We have to show that this is bimultiplicative, provided that  $(f\times f)^\ast{\cal M}$ is. Since it is symmetric, it suffices to show that $\cal M$ is a multiplicative sheaf when we view $V\times_SV$ as a group $V$-scheme via ${\rm pr}_1$. 
Consider the commutative diagram 
\begin{equation*}
\xymatrix{
U\times_SU\ar[d]_-{{\rm pr}_1}\ar[r]^-{{\rm id}\times f}&U\times_SV\ar[d]_-{{\rm pr}_1}\ar[r]^-{f\times{\rm id}}&V\times_SV\ar[d]_-{{\rm pr}_1}\\
U\ar[r]^-{\rm id}&U\ar[r]^-f&V. 
}
\end{equation*}
 By the assumption, $(f\times f)^\ast{\cal M}$ is a multiplicative sheaf when $U\times_SU$ is viewed as a group $U$-scheme via ${\rm pr}_1$. By the assertion $1$ applied to $U\times_SU\to U\times_SV$, we know that $(f\times{\rm id})^\ast{\cal M}$ is multiplicative. 
 Then Corollary \ref{fedes}.2 implies that $\cal M$ is multiplicative. The assertion follows. 
 \qed}

 \begin{cor}\label{desunip}
 Let $f\colon U\to V$ be a surjective morphism of perfect unipotent $k$-group schemes. Let $\cal L$ be an invertible $\Lambda$-sheaf on $V$. Then the following holds true. 
 \begin{enumerate}
 \item If $f^\ast\cal L$ is multiplicative, then $\cal L$ is multiplicative. 
 \item If $f^\ast\cal L$ is quadratic, then $\cal L$ is quadratic. 
 \end{enumerate}
 \end{cor}
 \proof{
 Let $f_0\colon U_0\to V_0$ be a model for $f$, as is constructed in Lemma \ref{redalg}. 
  As \'etale topoi do not change after the perfection, it suffices to show the similar statements for $f_0$. 
  Let $K_{\rm red}$ be the reduced closed subscheme of $U_0$ with the same underlying space as $\mathop{\rm ker}f_0$: this is a smooth closed group subscheme as $k$ is perfect. Then 
 the map $f_0$ is decomposed as follows: 
 \begin{equation*}
 U_0\to U_0/K_{\rm red}\to V_0, 
 \end{equation*}
 where the first map is smooth surjective and the second one is finite surjective radiciel. Therefore, it suffices to show the statements when $f_0$ has either  of the two properties. If $f$ is smooth surjective, then it is a consequence of Lemma \ref{desmq}. If $f_0$ is radiciel, then it follows as the \'etale topoi of $U_0,V_0$ are equivalent. 
 \qed}

\section{Computations of $H^i_c(U_{\bar{k}},{\cal Q})$}

Let $\bar{k}$ be an algebraic closure of $k$. 
Let $U$ be a perfect unipotent commutative group $k$-scheme and $\cal Q$ be a quadratic $\Lambda$-sheaf on $U$. In this section, we prove various properties of the cohomology groups of $\cal Q$ by using the results 
established in the previous sections and in the appendix. 
Let ${\cal B}_{\cal Q}$ be the bimultiplicative $\Lambda$-sheaf on $U\times_kU$ associated with $\cal Q$ (Definition \ref{quaddef}.1). 

\subsection{The rank of $H^i_c(U_{\bar{k}},{\cal Q})$: isogeneous case}\label{isogcase}
In this and next subsections, we compute the rank of $H^i_c(U_{\bar{k}},{\cal Q})$; notably, it turns out that the cohomology is concentrated in at most one degree. 

In this subsection, we treat the isogeneous case. 
Suppose that $\cal Q$ is isogeneous. 
To state Theorem \ref{HD},  we need to recall the following result. 
\begin{pr}(\cite[Corollary 6.9]{AD})\label{even}
Let $\cal B$ be a  bimultiplicative $\Lambda$-sheaf on $U\times_kU$ which is symmetric and isogeneous. Then its degree is of the form $p^{2r}$ for some integer $r\geq0$. 
\end{pr}
\proof{
See loc.~cit.~for the proof. Alternatively, one can use Corollary \ref{nondegdes} to prove the statement for those ${\cal B}$ that are of the form ${\cal B}_{\cal Q}$ for an isogeneous quadratic $\cal Q$. 
\qed}
 
\begin{thm}\label{HD}
Suppose that $\cal Q$ is isogeneous of degree $p^{2r}$. Then 
  the \'etale cohomology group $H^i_c(U_{\bar{k}},{\cal Q})$ vanishes unless $i=\dim U$. When $i=\dim U$, it is a finite free $\Lambda$-module of rank $p^r$. 
 %Suppose that $\Lambda=\overline{\mathbb{Q}_\ell}$ and that ${\cal Q}$ is an isogeneous quadratic $\Ql$-sheaf. Then the character ${\rm Gal}(\bar{k}/k)\to\Ql^\times$ corresponding to $({\rm det}H^{{\rm dim}U}_c(U_{\bar{k}},{\cal Q}))^2\otimes\Ql({\rm dim}U)$ has finite image. 
\end{thm}

%\begin{rmk} In proving the theorem, we may replace $k$ with arbitrary finite extension $k'/k$. Indeed, for the statement $1.$, this is obvious. For $2.$, let $\rho\colon {\rm Gal}(\bar{k}/k)\to\Ql^\times$ be a character and suppose that the image of ${\rm Gal}(\bar{k}/k')$ is finite. Then the image $\rho({\rm Gal}(\bar{k}/k))$ is also finite since the multiplication map ${\rm Gal}(\bar{k}/k)^{ab}\xrightarrow{[k':k]}{\rm Gal}(\bar{k}/k)^{ab}$ factors over the transfer map ${\rm Gal}(\bar{k}/k)^{ab}\to{\rm Gal}(\bar{k}/k')^{ab}$. In the sequel, we freely use this remark. 
%\end{rmk}

The aim of this subsection is to prove this theorem. 
We start with 
the following lemma. 

\begin{lm}\label{twist}
Let $\cal Q$ be an isogeneous quadratic $\Lambda$-sheaf on $U$ and let ${\cal L}$ be a multiplicative $\Lambda$-sheaf on $U$. Then 
${\cal Q}':={\cal Q}\otimes{\cal L}$ is again isogeneous quadratic. In the derived category of $\Lambda$-modules, $R\Gamma_c(U_{\bar{k}},{\cal Q}')$ is isomorphic to 
$R\Gamma_c(U_{\bar{k}},{\cal Q})$. 
\end{lm}
\proof{We may assume that $k=\bar{k}$. Let ${\cal B}_{\cal Q}$ denote the associated bimultiplicative $\Lambda$-sheaf and let $l_{{\cal Q}}\colon U\to U^\ast$ be the corresponding map. Consider  $\cal L$ as a $k$-rational point of $U^\ast$. Since $l_{{\cal Q}}$ is finite \'etale, we may find a  point $x\in U(k)$ mapping to $\cal L$. In other words, one can write ${\cal L}={\cal B}_{\cal Q}(-,x)$. Then we have 
\begin{equation*}
{\cal Q}'={\cal Q}\otimes{\cal B}_{\cal Q}(-,x)\cong {\cal Q}(-+x)\otimes\pi^\ast{\cal Q}(x)^{-1}. 
\end{equation*}
Here ${\cal Q}(-+x)$ denotes the pullback of $\cal Q$ via the translation map $U\to U,y\mapsto y+x$ and $\pi\colon U\to {\rm Spec}(k)$ is the structure map. Theorefore, we have  
\begin{equation*}
R\Gamma_c(U,{\cal Q}')\cong R\Gamma_c(U,{\cal Q}(-+x))\otimes{\cal Q}(x)^{-1}\cong R\Gamma_c(U,{\cal Q})\otimes{\cal Q}(x)^{-1}. 
\end{equation*}
Then the assertion follows. 
\qed
}

We prove Theorem \ref{HD} by induction on the dimension of $U$. First we prove the case when $U=\G_a^{\rm perf}$ by using the classification given in Lemma \ref{clB}. In the proof below, we identify \'etale sheaves on $\G_{a,k}^{\rm perf}$ and those on $\G_{a,k}$ as their \'etale topoi are canonically identified. 
Fix a non-trivial character $\psi\colon \F_p\to \Lambda^\times$. 
\begin{lm}\label{Gacase}
Theorem \ref{HD} holds true when $U=\G_{a,k}^{\rm perf}$. 
\end{lm}
\proof{
We may assume that $k=\bar{k}$. 
 Let $\cal Q$ be an isogeneous quadratic sheaf on $\G_{a,k}^{\rm perf}$. By Lemma \ref{clB}, ${\cal B}_{\cal Q}$ is isomorphic to ${\cal L}_\psi(f(x)y)$ for $f(X)=a_0X+\sum_{i\geq1}(a_iX^{p^i}+a_i^{p^{-i}}X^{p^{-i}})$. Since ${\cal B}_{\cal Q}$ is assumed to be isogeneous, $a_i$ is non-zero for some integer $i$. Let $n$ be the maximal integer with this property. Then one can rewrite $f(X)=a_0X+(a_1X^p+a_1^{p^{-1}}X^{p^{-1}})+\cdots+(a_nX^{p^n}+a_n^{p^{-n}}X^{p^{-n}})$. Our strategy to prove the statement is as follows. We will find another quadratic sheaf ${\cal Q}'$ such that ${\cal B}_{{\cal Q}'}\cong{\cal B}_{\cal Q}$. Moreover, ${\cal Q}'$ will be given in an explicit way. From its   explicit description, we are able to compute $R\Gamma_c(\G_{a,\bar{k}}^{\rm perf},{\cal Q}')$. Then the assertion for $\cal Q$ is reduced to that for ${\cal Q}'$ by Lemma \ref{twist}.

 We prove the statement by division into cases 
according to the parity of $p$. First suppose that $p$ is odd. Then we set ${\cal Q}':={\cal L}_\psi(\frac{1}{2}f(x)x)$. Using ${\cal L}_\psi(x^{p^{-i}+1})\cong {\cal L}_\psi(x^{p^{i}+1})$, we find that ${\cal Q}'\cong {\cal L}_\psi(\frac{1}{2}a_0x^2)\otimes(\bigotimes_{i=1}^n{\cal L}_\psi(a_ix^{p^i+1}))$. Note that the Swan conductor of the locally constant sheaf ${\cal L}_\psi(a_ix^{p^i+1})$ on $\G_{a,k}$ at the infinity is equal to $p^i+1$ if $a_i\neq0$ and $0$ if $a_i=0$. Therefore, that of ${\cal Q}'$ is  equal to $p^n+1$. Since ${\cal Q}'$ and ${\cal Q}'^{-1}(1)$ ($(-)$ means Tate twist) are non-trivial, the cohomology groups $H^0_c(\G_{a,k},{\cal Q}')$ and 
$H^2_c(\G_{a,k},{\cal Q}')\cong H^0(\G_{a,k},{\cal Q}'^{-1}(1))^\vee$ vanish. 
Then  the Grothendieck--Ogg--Shafarevich formula implies that $H^1_c(\G_{a,k},{\cal Q}')$ has rank $p^n$, as claimed. 

We treat the remaining case where $p=2$. 
Let $W_2(-)$ denote the $2$-typical Witt ring of length $2$. 
Let $\psi'\colon W_2(\F_2)\to\Lambda^\times$ be a character with the property that the composite map  $\F_2\xrightarrow{V}W_2(\F_2)\to\Lambda^\times$ is equal to $\psi$. 
Let $W_2$ denote the affine ring scheme over $k$ which represents the presheaf $Sch/k\ni S\mapsto W_2(\Gamma(S,{\cal O}_S))$. The character $\psi'$ and the \'etale covering $W_2\to W_2$ defined by $(x,y)\mapsto (x^2,y^2)-(x,y)$ induces a locally constant sheaf on $W_2$. 
We write ${\cal L}_{\psi'}$ for this sheaf.
 Let ${\cal L}_{\psi'}(\sqrt{a_0}x,0)$ be its pullback via the map $\G_{a,k}\to W_2, x\mapsto (\sqrt{a_0}x,0)$. We consider the quadratic sheaf ${\cal Q}':={\cal L}_{\psi'}(\sqrt{a_0}x,0)\otimes(\otimes_{i=1}^n{\cal L}_\psi(a_ix^{p^i+1}))$. It is straightforward to check that ${\cal B}_{{\cal Q}'}={\cal B}_{\cal Q}$. If $n\geq1$, then the Swan conductor of ${\cal Q}'$ is equal to $p^n+1$. If $n=0$, then ${\cal Q}'$ is isomorphic to ${\cal L}_{\psi'}(\sqrt{a_0}x,0)$, which has  Swan conductor $2$. Then the Grothendieck--Ogg--Shafarevich formula applies. 
\qed}

\vspace{5mm}

(Proof of Theorem \ref{HD})

\proof{We may assume that $k=\bar{k}$. 
We argue by induction on $\dim U$. When $
\dim U=1$, $U$ is isomorphic to $\G_{a,k}^{\rm perf}$, for which the assertion is proved in Lemma \ref{Gacase}. 

Suppose that $\dim U>1$. Take and fix a closed subgroup $V\subset U$ which is isomorphic to $\G_{a,k}^{\rm perf}$. Let $V^\perp$ be the kernel of the map $U\xrightarrow{l_{{\cal B}_{\cal Q}}}U^\ast\to V^\ast$. 
By Lemma \ref{a}, this is a commutative group scheme over $k$ with finitely many connected components. Moreover, its connected component of the neutral element, for which we write $V^{\perp0}$, 
is perfect unipotent commutative. 
Let ${\cal Q}_V$ be the restriction of $\cal Q$ to $V$.
This is a quadratic sheaf on $V$. 
 We have two cases: ${\cal Q}_V$ is isogeneous or it is not. Note that, in the latter case,  ${\cal B}_{\cal Q}$ is trivial, which means that ${\cal Q}_V$ is a multiplicative $\Lambda$-sheaf. 

First suppose that ${\cal Q}_V$ is isogeneous. Then the intersection $V\cap V^{\perp0}$ is finite and the map $\pi\colon V\times V^{\perp0}\to U$ given by $(x,y)\mapsto x+y$ is finite \'etale by Lemma \ref{morper}.2. 
Thus $\pi^\ast{\cal Q}$ is isogeneous with degree equal to $\deg{\cal Q}\cdot(\deg\pi)^2$ by Lemma \ref{pullQ}.2. Therefore, we know that ${\cal Q}|_{V^{\perp0}}$ is an isogeneous quadratic sheaf by 
noting that $\pi^\ast{\cal Q}$ is isomorphic to ${\cal Q}_V\boxtimes{\cal Q}|_{V^{\perp0}}$ since ${\cal B}_{\cal Q}|_{V\times V^{\perp0}}$ is trivial. Moreover, we have the following equality on the degrees: 
\[\deg{\cal Q}\cdot(\deg\pi)^2=\deg{\cal Q}_V\cdot\deg{\cal Q}_{V^{\perp0}}. 
\]

Let $A$ be the kernel of $\pi$, which we identify with the finite abelian group $A(k)$. For $\chi\in A^\vee$, let ${\cal L}_\chi$ be the invertible \'etale $\Lambda$-sheaf on $U$ associated with $\chi$ and the \'etale covering $\pi$: this is a multiplicative $\Lambda$-sheaf. We have 
\begin{align*}
R\Gamma_c(V,{\cal Q}_V)\otimes R\Gamma_c(V^{\perp0},{\cal Q}|_{V^{\perp0}})&\cong R\Gamma_c(V\times V^{\perp0}, \pi^\ast {\cal Q})\\
&\cong \bigoplus_{\chi\in A^\vee}R\Gamma_c(U,{\cal Q}\otimes{\cal L}_\chi). 
\end{align*}
Then the assertion follows from Lemma \ref{twist}, \ref{Gacase}, and the induction hypothesis applied to $(V^{\perp0},{\cal Q}|_{V^{\perp0}})$.

Next suppose that 
${\cal Q}_V$ is multiplicative. We claim that there exists a multiplicative $\Lambda$-sheaf $\cal L$ on $U$ such that ${\cal L}|_V\cong {\cal Q}_V$. Indeed, by Proposition \ref{dualp}.2, the sequence 
\[0\to(U/V)^\ast\to U^\ast\to V^\ast\to0
\]
is exact. Then, by Lemma \ref{section}, the map $U^\ast\to V^\ast$ has a section, which gives an extension $\cal L$. 
By replacing $\cal Q$ with ${\cal Q}\otimes{\cal L}^{-1}$ (this is allowed by Lemma \ref{twist}), we may assume that ${\cal Q}_V$ is trivial. 

Suppose that ${\cal Q}_V$ is trivial.  Let $f\colon U\to U/V$ be the quotient map. We compute the pushforward $Rf_!\cal Q$. For a point $a\in U(k)$, let ${\cal Q}_{V,a}$ be the pullback of $\cal Q$ via $V\to U, x\mapsto x+a$. This is isomorphic to ${\cal B}_{\cal Q}(-,a)\otimes{\cal Q}(a)$ since ${\cal Q}_V$ is trivial. As $k$ is algebraically closed,  ${\cal Q}(a)$ is also trivial. Therefore, ${\cal Q}_{V,a}$ is a multiplicative $\Lambda$-sheaf which is trivial if and only if $a\in V^\perp$. Since the cohomology of a non-trivial multiplicative $\Lambda$-sheaf on $\G_{a,k}^{\rm perf}$ is acyclic, $Rf_!{\cal Q}$ is supported on $V^\perp/V$. Let $\pi_0(V^\perp)$ be the set of connected components and, 
for $\alpha\in \pi_0(V^\perp)$, let $V^{\perp\alpha}$ be the connected component corresponding to $\alpha$. We have 
\begin{equation}\label{QRf}
R\Gamma_c(U,{\cal Q})\cong \bigoplus_{\alpha\in\pi_0(V^\perp)}R\Gamma_c(V^{\perp\alpha}/V,Rf_!{\cal Q}). 
\end{equation}
For each $\alpha$, choose an element $a\in V^{\perp\alpha}(k)$ and let $t_a\colon V^{\perp0}\to 
V^{\perp\alpha}$ be the map given by the translation by $a$. Then  we have 
\[ t_a^\ast{\cal Q}|_{V^{\perp\alpha}}\cong {\cal Q}|_{V^{\perp0}}\otimes{\cal B}_{\cal Q}(-,a). 
\]
Note that this is a quadratic $\Lambda$-sheaf on $V^{\perp0}$. We claim that this sheaf descends to a sheaf on $V^{\perp0}/V$. To show this,  consider the fiber product 
\begin{equation*}
\xymatrix{
V^{\perp0}\times V\ar[d]_-{{\rm pr}_1}\ar[r]^-{m}& V^{\perp0}\ar[d]\\
V^{\perp0}\ar[r]&V^{\perp0}/V.
}
\end{equation*}
Here $m$ is the map $(x,y)\mapsto x+y$. 
We note that $m^\ast t_a^\ast{\cal Q}|_{V^{\perp\alpha}}$ is isomorphic to ${\rm pr}_1^\ast t_a^\ast{\cal Q}|_{V^{\perp\alpha}}$. Indeed, $m^\ast{\cal Q}|_{V^{\perp0}}$ is isomorphic to ${\rm pr}_1^\ast{\cal Q}|_{V^{\perp0}}$ 
since ${\cal B}_{\cal Q}|_{V^{\perp0}\times V}$ and ${\cal Q}_V$ are trivial. Moreover, $m^\ast{\cal B}_{\cal Q}(-,a)$ is isomorphic to ${\rm pr}_1^\ast{\cal B}_{\cal Q}(-,a)$ as $a\in V^\perp$. Then  Lemma \ref{des} shows that there exists an invertible \'etale $\Lambda$-sheaf ${\cal Q}'_a$ on $V^{\perp0}/V$ that is pulled back to $t^\ast_a{\cal Q}|_{V^{\perp\alpha}}$. By Corollary \ref{desunip}, we moreover know that this is a quadratic $\Lambda$-sheaf. 

The isomorphism (\ref{QRf}) can be rewritten as follows: 
\[
 R\Gamma_c(U,{\cal Q})\cong \bigoplus_{\alpha\in\pi_0(V^\perp)}R\Gamma_c(V^{\perp0}/V,{\cal Q}'_a)(-1)[-2].
 \]
 To complete the proof, it suffices to show that ${\cal Q}'_a$ is isogeneous with degree equal to $\deg{\cal Q}/|\pi_0(V^\perp)|$. This assertion is reduced to showing the following three equalities. Let $B$ be the kernel of $l_{\cal Q}$. 
 \begin{enumerate}
 \item[(1)] $|\pi_0(V^\perp)|=|B\cap V|$. 
 \item[(2)] $|\pi_0(V^\perp)|=|B|/|B\cap V^{\perp0}|$. 
 \item[(3)] $\deg{\cal Q}_a'=|B\cap V^{\perp0}|/|B\cap V|$. 
 \end{enumerate}
 We prove (1). The map $U\to V^\ast$ induced by $l_{\cal Q}$ factors over a map $U/V^{\perp0}\to V^\ast$. This is finite \'etale and its degree is equal to $|\pi_0(V^\perp)|$. By Theorem \ref{dualp}.3, this degree is equal to that of $V\to (U/V^{\perp0})^\ast$. 
 Since the kernel of the latter map is $B\cap V$, we get the equality (1). 
 
 We prove (2). We have a cartesian square 
 \begin{equation*}
 \xymatrix{
 U\ar[r]^-{l_{\cal Q}}&U^\ast\\
 V^\perp\ar[u]\ar[r]&(U/V)^\ast,\ar[u]
 }
 \end{equation*}
 where the bottom arrow is surjective. As $(U/V)^\ast$ is connected, the map $V^{\perp0}\to(U/V)^\ast$ is also surjective. Thus we have a commutative diagram 
 \begin{equation}\label{experp}
 \xymatrix{
 0\ar[r]&B\ar[r]&V^\perp\ar[r]&(U/V)^\ast\ar[r]&0\\
 0\ar[r]&B\cap V^{\perp0}\ar[r]\ar[u]&V^{\perp0}\ar[u]\ar[r]&(U/V)^\ast\ar[r]\ar[u]^-{{\rm id}}&0
 }
 \end{equation}
 where the horizontal sequences are exact and the vertical arrows are the inclusions. Then, by the snake lemma, we get the equality (2). 
 
 We prove (3). As the quotient map $V^{\perp0}\to V^{\perp0}/V$ has geometrically connected fibers, $\deg{\cal Q}'_a$ is equal to the number of connected components of the kernel of the composite map $V^{\perp0}\to V^{\perp0}/V\xrightarrow{l_{{\cal Q}'_a}}
 (V^{\perp0}/V)^\ast$. This map has another factorization $V^{\perp0}\to (U/V)^\ast\to(V^{\perp0}/V)^\ast$, where the former map is the one appearing in (\ref{experp}). The kernel of $(U/V)^\ast\to(V^{\perp0}/V)^\ast$ is equal to $(U/V^{\perp0})^\ast$ and the inverse image of $(U/V^{\perp0})^\ast$ to $V^{\perp0}$ is $V+(B\cap V^{\perp0})$. 
 Then the equality (3) follows by using $(V+(B\cap V^{\perp0}))/V\cong (B\cap V^{\perp0})/(B\cap V).$
 
 The proof is completed. 
 \qed}

\subsection{The rank of $H^i_c(U_{\bar{k}},{\cal Q})$: general case}\label{genrank}

Now suppose that $\cal Q$ is not necessarily isogeneous. The purpose of this subsection is to show that the computation of $H^i_c(U_{\bar{k}},{\cal Q})$ can be reduced to the isogeneous case. 

Let ${\cal B}_{\cal Q}$ be the bimultiplicative $\Lambda$-sheaf on $U\times_kU$ associated to $\cal Q$ and let $l_{\cal Q}\colon U\to U^\ast$ be the corresponding group homomorphism. Let $K^\circ$ be the neutral connected component of $\mathop{\rm ker}l_{\cal Q}$, which is perfect unipotent by Lemma \ref{a}. Then the restriction ${\cal Q}|_{K^\circ}$ is a multiplicative $\Lambda$-sheaf on $K^\circ$. 
 Let $\pi\colon U\to U/K^\circ$ be the quotient map. 
 There are two possible cases. 
 
 \vspace{2mm}
 
 $(1)$ the restriction ${\cal Q}|_{K^\circ}$ is non-trivial. 
 
 \vspace{2mm}
 
 In this case, the pushforward $R\pi_!{\cal Q}$ is acyclic. Indeed, for any $\bar{k}$-rational point $a\in (U/K^\circ)(\bar{k})$, the stalk $(R\pi_!{\cal Q})_a$ is isomorphic to $R\Gamma_c(K^\circ_{\bar{k}},{\cal Q}|_{K^\circ})\otimes{\cal Q}(a)$, which is acyclic by Lemma \ref{vancoh}. Consequently, we have $R\Gamma_c(U_{\bar{k}},{\cal Q})\cong 0$ in this case. 
 
 \vspace{2mm}
 
 $(2)$ the restriction ${\cal Q}|_{K^\circ}$ is trivial. 
 
 \vspace{2mm}
 \begin{lm}
 Suppose that ${\cal Q}|_{K^\circ}$ is trivial. Then there exists an isogeneous quadratic $\Lambda$-sheaf $\overline{{\cal Q}}$ on $U/K^\circ$ such that $\pi^\ast 
 \overline{{\cal Q}}$ is isomorphic to $\cal Q$. 
 \end{lm}
 \proof{
 Consider the cartesian square
 \begin{equation*}
 \xymatrix{
 U\times_kK^\circ\ar[d]_-{{\rm pr}_1}\ar[r]^-m&U\ar[d]^-\pi\\
 U\ar[r]_-\pi&U/K^\circ,
 }
 \end{equation*}
 where $m$ is given by $(x,y)\mapsto x+y$. 
 The triviality of ${\cal Q}|_{K^\circ}$ implies that $m^\ast {\cal Q}$ is isomorphic to ${\rm pr}_1^\ast{\cal Q}$. By Lemma \ref{des}, there exists an invertible \'etale $\Lambda$-sheaf $\overline{{\cal Q}}$ whose pullback by $\pi$ is isomorphic to $\cal Q$. By Corollary \ref{desunip}.2, $\overline{{\cal Q}}$ is quadratic. 
 We have a commutative diagram 
 \begin{equation}\label{UKQ}
 \xymatrix{
 U\ar[d]_-{l_{\cal Q}}\ar[r]^-\pi&U/K^\circ\ar[d]^-{l_{\overline{{\cal Q}}}}\\
 U^\ast&(U/K^\circ)^\ast\ar[l]^-{\pi^\ast}. 
 }
 \end{equation}
 Since $\pi$ is surjective and $\pi^\ast$ is a closed immersion, $l_{\overline{{\cal Q}}}$ must be an isogeny by reason of the dimension of the kernel. 
 \qed}

We can summarize the results as follows. 
\begin{cor}\label{rankgen}
Let $\cal Q$ be a quadratic $\Lambda$-sheaf on $U$. Let $K$ be the kernel of $l_{\cal Q}$ and $K^\circ$ be its neutral connected component. 
\begin{enumerate}
\item Suppose that ${\cal Q}|_{K^\circ}$ is non-trivial. Then we have $H^i_c(U_{\bar{k}},{\cal Q})\cong0$ for all $i$. 
\item Suppose that ${\cal Q}|_{K^\circ}$ is trivial. Then $H^i_c(U_{\bar{k}},{\cal Q})$ is trivial except for $i=\dim(U/K^\circ)$. 
When $i=\dim(U/K^\circ)$, it is a free $\Lambda$-module with rank equal to the cardinality of $\pi_0(K)$. 
\end{enumerate}
\end{cor}
\proof{
The only non-trivial part is the last assertion in $2.$ 
Consider the diagram (\ref{UKQ}). Since $\pi^\ast$ is a closed immersion and $\pi$ has geometrically connected fibers, the cardinality of $\mathop{\rm ker}l_{\overline{\cal Q}}$ is equal to the number of connected components of $\mathop{\rm ker}l_{\cal Q}$, hence the assertion. 
\qed}

\section{Heisenberg group attached to $\cal Q$}
\label{secH}
In this section, we observe that one can naturally construct a certain finite Heisenberg group from a quadratic sheaf $\cal Q$ and prove various results on $\cal Q$ which can be deduced by considering the Heisenberg group.

We use the following notation throughout this section: $U$ is a perfect unipotent commutative group $k$-scheme and $\cal Q$ is a quadratic $\Lambda$-sheaf on $U$. Let $C\to U$ be the finite Galois covering that trivializes $\cal Q$, i.e., the character $\pi_1(U)^{ab}\to\Lambda^\times$ corresponding to $\cal Q$ factors over an injection ${\rm Aut}(C/U)\to\Lambda^\times$. 
Let $\bar{k}$ be an algebraic closure of $k$. 
We assume that $\cal Q$ is isogeneous unless otherwise stated. 

\subsection{Construction of Heisenberg group and its basic properties}\label{H}
In this subsection, we suppose that $\cal Q$ is isogeneous. 
Let $l_{\cal Q}\colon U\to U^\ast$ be the map associated with ${\cal B}_{\cal Q}$, which is finite \'etale. The aim of this subsection is to show that the composite map $C\to U\to U^\ast$ is a Galois covering at least after a finite base field extension and to investigate some basic  properties of its Galois group. We start with the following. 
\begin{lm}\label{Heis}
The map $C_{\bar{k}}\to U^\ast_{\bar{k}}$ is a Galois covering. 
\end{lm}
\proof{
We show that the cardinality of ${\rm Aut}(C_{\bar{k}}/U^\ast_{\bar{k}})$ is equal to the degree of $C\to U^\ast$. Take $a\in \mathop{\rm ker}l_{\cal Q}(\bar{k})$. For any point $x\in U(S)$ for $S\in Psch/\bar{k}$, we have an isomorphism 
\begin{equation*}
{\cal Q}(x+a)\cong{\cal Q}(x)\otimes{\cal Q}(a)\otimes{\cal B}_{\cal Q}(x,a). 
\end{equation*}
Note that ${\cal B}_{\cal Q}(x,a)$ is trivial since $a\in\mathop{\rm ker}l_{\cal Q}(\bar{k})$. Also ${\cal Q}(a)$ is trivial. Therefore, we have an isomorphism ${\cal Q}(x+a)\cong{\cal Q}(x)$, which implies that the automorphism $U_{\bar{k}}\to U_{\bar{k}},x\mapsto x+a$ can be extended to an automorphism of $C_{\bar{k}}$. 
\qed}

Put $A:={\rm Aut}(C_{\bar{k}}/U_{\bar{k}})(\cong {\rm Aut}(C/U))$, 
 $H:={\rm Aut}(C_{\bar{k}}/U^\ast_{\bar{k}})$, and $ K:={\rm Aut}(U_{\bar{k}}/U^\ast_{\bar{k}})=\mathop{\rm ker}l_{\cal Q}(\bar{k})$. Note that $A$ is a cyclic group of $p$-power order (Lemma \ref{order}.2). 
 We have a short exact sequence of groups 
\begin{equation*}
1\to A\to H\to K\to 1. 
\end{equation*}

\begin{lm}\label{centerA}
The following hold true. 
\begin{enumerate}
\item We have an inclusion $A\subset Z(H)$. Here $Z(H)$ denotes the center. 
\item The map $H\times H\to A, (x,y)\mapsto xyx^{-1}y^{-1}$ induces an alternating bilinear pairing $\langle-,-\rangle\colon K\times K\to A$. This pairing induces an isomorphism $K\to{\rm Hom}(K,A)$. 
\item We have an equality $A=Z(H)$. 
\end{enumerate}
\end{lm}
A finite group $G$ with the property that $G/Z(G)$ is abelian is known as a Heisenberg group. For properties of such groups which we use in the sequel, we refer to the book \cite{Bump}, Exercises 4.1.4--4.1.8. 
 \proof{
 We write the quotient map $\pi_1(U_{\bar{k}})^{ab}\to A$ as $\sigma\mapsto\bar{\sigma}$. 
 
$1$.  We prove that $A\subset Z(H)$. Take any $\alpha \in H$ and let $a\in K=\mathop{\rm ker}l_{\cal Q}$ be its image. Let $ \pi_1(U)^{ab}\to A\xrightarrow{\rho_{\cal Q}}\Lambda^\times$ be the character corresponding to $\cal Q$. Then the character corresponding to ${\cal Q}(x+a)$ is given by $\pi_1(U_{\bar{k}})^{ab}\ni \sigma\mapsto\rho_{\cal Q}(\alpha \bar{\sigma}\alpha^{-1})$. Since ${\cal Q}(x+a)$ is isomorphic to ${\cal Q}(x)$, this character is equal to $\rho_{\cal Q}$. As $\rho_{\cal Q}$ is injective, we must have $A\subset Z(H)$. 
  
 $2.$ It is straightforward to check that the map $H\times H\to A$ induces an alternating pairing $\langle-,-\rangle$. 
 Let $M:=\{x\in K\mid\forall y\in K, \langle x,y\rangle=1\}$. We show that this is trivial as follows. 
 Since $A$ is cyclic, the map $K/M\to {\rm Hom}(K/M,A)$ induced by the pairing is an isomorphism. Take a maximal totally isotropic subgroup $B\subset K/M$. 
 Since $\langle-,-\rangle$ is alternating, we have an equality $B=B^\perp$ where $B^\perp:=\{x\in K/M|\forall y\in B,\langle x,y\rangle=1\}$ (for this, see \cite[Ex.4.1.5]{Bump}), and consequently $|B|^2=|K/M|$. 
 
 Let $\widetilde{B}$ be the inverse image of $B$ via the quotient map $H\to K/M$ and let $V\to U^\ast_{\bar{k}}$ be the corresponding \'etale covering (i.e., ${\rm Aut}(C_{\bar{k}}/V)=\Tilde{B}$). Since $V$ sits between $U_{\bar{k}}\to U^\ast_{\bar{k}}$, $V$ naturally carries a structure of perfect unipotent group $\bar{k}$-scheme and $U_{\bar{k}}\to V_{\bar{k}}$ is a group homomorphism. Since $\widetilde{B}$ is an abelian group, the character $\rho_{\cal Q}\colon A\to \Lambda^\times$ can be extended to a character $\widetilde{B}\to\Lambda^\times$. Let ${\cal Q}'$ be the corresponding invertible \'etale sheaf on $V$. Since ${\cal Q}'$ is restricted to $\cal Q$, Corollary \ref{desunip}.2 implies that ${\cal Q}'$ is quadratic. Let $l_{{\cal Q}'}\colon V\to V^\ast$ be the map associated with ${\cal B}_{{\cal Q}'}$. This fits into the following commutative diagram
 \begin{equation*}
 \xymatrix{
 U_{\bar{k}}\ar[r]\ar[d]_-{l_{\cal Q}}&V\ar[d]^-{l_{{\cal Q}'}}\\ U^\ast_{\bar{k}}&V^\ast.\ar[l]. 
 }
 \end{equation*}
 From this, we have equalities $\deg l_{\cal Q}=\deg l_{{\cal Q}'}\cdot\deg(U_{\bar{k}}/V)^2$. This can be rewritten as 
 $|K|=\deg l_{{\cal Q}'}\cdot(|M||B|)^2$. Since $|B|^2=|K/M|$, we have 
 \begin{equation}\label{deg1}
 |K|=\deg l_{{\cal Q}'}\cdot|M|^2|K/M|=\deg l_{{\cal Q}'}\cdot|M||K|. 
 \end{equation}
 This implies that $|M|=1$, as desired. 
 
 $3.$ Let $x\in Z(H)$ and let $\bar{x}$ be its class in $H/A=K$. From the definition of the pairing $\langle -,-\rangle$, we have $\langle \bar{x},y\rangle=1$ for any $y\in K$. As $\langle -,-\rangle$ is a perfect pairing, this implies that $\bar{x}=1$, i.e., $x\in A$.  
 \qed}

 From the proof of Lemma \ref{centerA}, we have the following surprising result. 
 \begin{cor}\label{nondegdes}
 Let $\widetilde{B}$ be a maximal abelian subgroup of $H$. Suppose that every automorphism belonging to $\widetilde{B}$ is defined over $k$. Then the quotient scheme $V:=C/\widetilde{B}$ has a unique structure of perfect unipotent commutative group $k$-scheme such that $U\to V$ is an \'etale isogeny. Moreover, the quadratic sheaf $\cal Q$ can be obtained as the pullback of a non-degenerate quadratic $\Lambda$-sheaf on $V$.  
 \end{cor}
 \proof{
 The group $\widetilde{B}/A$ is naturally identified with a finite subgroup of $U(k)$. Therefore, $V$ carries a structure of perfect unipotent group. The uniqueness is clear as $U\to V$ is \'etale and surjective. By the same way as in the proof of Lemma \ref{centerA}.2, $V$ has a quadratic $\Lambda$-sheaf ${\cal Q}'$ whose pullback to $U$ is $\cal Q$. By (\ref{deg1}), we have $\deg l_{{\cal Q}'}=1$, which implies the non-degeneracy. 
 \qed}
 
 We remark the following well-known property of Heisenberg groups. 
 \begin{rmk}
 From the non-degeneracy of $\langle-,-\rangle$, we know that the assignement $\widetilde{B }\mapsto \widetilde{B}/A$ gives a bijection between the set of abelian subgroups of $H$ and that of totally isotropic subgroups of $H/A$. In particular, the maximal ones are in one to one correspondence. 
 \end{rmk}
 
 Let $\pi\colon C_{\bar{k}}\to U_{\bar{k}}$ denote the quotient map. This is equivariant with respect to the actions of $H$. The pushforward $\pi_\ast\Lambda$ has the following direct sum decomposition
 \[
\pi_\ast\Lambda\cong\bigoplus_{\psi\in A^\vee}{\cal Q}_\psi, 
 \]
 where  ${\cal Q}_\psi$ denotes the $\psi$-isotipic part. If $\psi$ is the injection $A\to\Lambda^\times$ corresponding to $\cal Q$, then ${\cal Q}_\psi$ is nothing but $\cal Q$. 
  As $A$ is the center of $H$, the action of $H$ on $\pi_\ast\Lambda$ respects this decomposition. In particular, $H$ acts on the cohomology group of $\cal Q$. 
 \begin{pr}\label{irrH}
 Further suppose that $\Lambda=\Ql$. Let $\psi\in A^\vee$ be the character given by $\cal Q$. 
 Then the representation of $H$ given by $H_c^d(U_{\bar{k}},{\cal Q})$, where $d=\dim U$, is a unique irreducible representation with central character $\psi$. 
 \end{pr}
\proof{ The Stone--Von Neumann theorem (\cite[Ex. 4.1.8]{Bump}) implies that there exists a unique irreducible representation $\pi_\psi$ of $H$ whose central character is $\psi$. Moreover, its dimension is known to be equal to $\sqrt{|H/A|}=\sqrt{|K|}$.  On the other hand,  $H_c^d(U_{\bar{k}},{\cal Q})$ is an $H$-representation on which $A$ acts as the scalar multiplication given by $\psi$. Thus $H_c^d(U_{\bar{k}},{\cal Q})$ is isomorphic to the direct
sum of finitely many copies of $\pi_\psi$. Since we have $\dim H_c^d(U_{\bar{k}},{\cal Q})=\sqrt{|K|}$ by Theorem \ref{HD}, the multiplicity should be $1$, i.e., 
the cohomology is isomorphic to $\pi_\psi$. 
 \qed}
 
 The absolute Galois group ${\rm Gal}(\bar{k}/k)$ acts on $H={\rm Aut}(C_{\bar{k}}/U^\ast_{\bar{k}})$ via conjugation. This action is trivial if and only if $C\to U^\ast$ is a Galois covering over $k$. 
 The irreducibility of the $H$-representation $H^d_c(U_{\bar{k}},{\cal Q})$ gives us the following criterion for whether $C\to U^\ast$ is Galois over $k$ or not. 
\begin{pr}\label{GandH}
Further suppose that $\Lambda=\Ql$. 
The following are equivalent. 
\begin{enumerate}
\item The \'etale covering $C\to U^\ast$ is a Galois covering. In other words, all automorphisms in $H$ can be defined over $k$. 
\item Any element in ${\rm Gal}(\bar{k}/k)$ acts as a scalar multiplication on $H^d_c(U_{\bar{k}},{\cal Q})$ where $d=\mathop{\rm dim}U$. 
\end{enumerate}
\end{pr}
\proof{
Suppose that $C\to U^\ast$ is Galois over $k$. Then the action of ${\rm Gal}(\bar{k}/k)$ commutes with that of $H$. As $H^d_c(U_{\bar{k}},{\cal Q})$ is an irreducible representation of $H$, this implies $2.$ 

 We prove the inverse. Let $\sigma\in{\rm Gal}(\bar{k}/k)$. We have to show that the map $H\to H,x\mapsto \sigma x\sigma^{-1}$ is equal to the identity map provided that $\sigma$ acts  as a scalar multiplication on $H^d_c(U_{\bar{k}},{\cal Q})$. To do so, it suffices to show that the map $H\to {\rm GL}(H^d_c(U_{\bar{k}},{\cal Q}))$ is injective, which follows from the following lemma. 
\qed}

\begin{lm}
Let $H$ be a finite group with the following properties: the center $A$ is cyclic and $H/A$ is abelian. 

Let $\psi\colon A\to \Ql^\times$ be an injective character. Let $\pi_\psi$ be the unique irreducible representation of $H$ whose central character is $\psi$ 
and  $V$ be the underlying vector space of $\pi_\psi$. Then the map $\rho_{\pi_\psi}\colon H\to {\rm GL}(V)$ given by $\pi_\psi$ is injective. 
\end{lm}
\proof{
Let $a\in H\setminus\{1\}$. We show that $\rho_{\pi_\psi}(a)\neq1$. If $a\in A$, then $\rho_{\pi_\psi}(a)$ is the multiplication by $\psi(a)$, which is not $1$ by the injectivity of $\psi$. Suppose that $a\notin A.$ Take a maximal abelian subgroup $M$ of $H$ that contains $a$. We have a short exact sequence of the character groups 
\begin{equation*}
0\to(M/A)^\vee\to M^\vee\to A^\vee\to0. 
\end{equation*}
We claim that there exists $\psi'\in M^\vee$ such that $\psi'|_A=\psi$ and $\psi'(a)\neq1$. Indeed, 
take one $\psi'\in M^\vee$ such that $\psi'|_A=\psi$. If $\psi'(a)\neq1$, then we are done. Suppose $\psi'(a)=1$. Choose $\eta\in(M/A)^\vee$ such that $\eta(\bar{a})\neq1$, where $\bar{a}$ is the image of $a$ in $M/A$. This is possible since $\bar{a}\neq0$ in $M/A$. Then $\psi'\eta$ satisfies this property. 

Since the restriction $\pi_\psi|_M$ is isomorphic to $\bigoplus_{\psi'}\psi'$ where $\psi'$ runs through the characters $\psi'\in M^\vee$ such that $\psi'|_A=\psi$, $a$ acts non-trivially on $V$. 
\qed}

\subsection{Frobenius eigenvalues of $H^d_c(U_{\bar{k}},{\cal Q})$}\label{H2}

Suppose that $k$ is a finite field $\F_q$ with $q$ elements. In this subsection, we prove a result on the Frobenius  eigenvalues, which can be regarded as a cohomological interpretation of Proposition \ref{GS}. 

\begin{pr}\label{FEV}
Let $U$ be a perfect unipotent commutative group $\F_q$-scheme and $\cal Q$ be an isogeneous quadratic $\Ql$-sheaf on it. Let $d=\dim U$. Then any  eigenvalue of the geometric Frobenius ${\rm Fr}_q$ on $H^d_c(U_{\overline{\F_q}},{\cal Q})$ is of the form $\zeta\cdot\sqrt{q^d}$ where $\zeta$ is a root of unity. 
\end{pr}

Combining this proposition and results established above, we have the following version of Hasse--Davenport relation. To state it, we fix a notation for the function--sheaf dictionary: for a scheme $X$ of finite type over $\F_q$ and a $\Ql$-sheaf $\cal F$ on it, we write the function $X(\F_{q^n})\to\Ql$ given by $x\mapsto {\rm Tr}({\rm Fr}_{q^n},{\cal F}_{\bar{x}})$, where ${\rm Fr}_{q^n}$ denotes the $q^n$-th geometric Frobenius, as $t_{{\cal F},q^n}$. 

\begin{cor}\label{genHD}
We follow the notation and assumptions made  in Proposition \ref{FEV}. Let $p^{2r}$ denote the degree of the map $\l_{\cal Q}$. 

Let $C\to U$ be the Galois covering that trivializes $\cal Q$. We assume that any element in $H={\rm Aut}(C_{\overline{\F_q}}/U^\ast_{\overline{\F_q}})$ is defined over $\F_q$. Then there exists an element $\tau\in\Ql^\times$ which satisfies the following relation: for any integer $n\geq1$, 
\[\sum_{x\in U(\F_{q^n})}t_{{\cal Q},q^n}(x)=(-1)^dp^r\tau^n. 
\]
Moreover, $\tau$ is of the form $\zeta\cdot\sqrt{q^d}$ where $\zeta$ is a root of unity. 
\end{cor}
Consequently, if $\cal Q$ is non-degenerate, then we have the following relation for any $n\geq1$: 
\[(-1)^d\sum_{x\in U(\F_{q^n})}t_{{\cal Q},q^n}(x)
=((-1)^d\sum_{x\in U(\F_{q})}t_{{\cal Q},q}(x))^n.\]
Indeed, in the non-degenerate case, the covering $C\to U^\ast$ is always Galois as $l_{\cal Q}$ is an isomorphism and $C\to U$ is Galois. 

\proof{
By Grothendieck's trace formula, we have 
\[\sum_{x\in U(\F_{q^n})}t_{{\cal Q},q^n}(x)=\sum_{i\geq0}(-1)^i{\rm Tr}({\rm Fr}_{q^n}, H^i_c(U_{\overline{\F_q}},{\cal Q})). 
\]
By Theorem \ref{HD}, $H^i_c(U_{\overline{\F_q}},{\cal Q})$ vanishes if $i\neq d$ and has dimension $p^r$ if $i=d$. Moreover, Proposition \ref{GandH} implies that ${\rm Fr}_q$ acts as a scalar multiplication by some $\tau\in\Ql^\times$. The last assertion follows from Proposition \ref{FEV}. 
\qed}

The last of this subsection is devoted to proving Proposition \ref{FEV}. We start with the case where $\cal Q$ is non-degenerate, for which case we can apply Proposition \ref{GS}. 

\begin{pr}\label{nondegcase}
Under the notation and assumption in Proposition \ref{FEV}, we further assume that $\cal Q$ is non-degenerate. 
\begin{enumerate}
\item The mapping $t_{{\cal Q},q}\colon U(\F_q)\to\Ql^\times$ defined by $x\mapsto {\rm Tr}({\rm Fr}_q, {\cal Q}_{\bar{x}})$ is a non-degenerate qudratic form in the sense of Definition \ref{quadcl}. 
\item Let $\tau$ be the eigenvalue of ${\rm Fr}_q$ on $H^d_c(U_{\overline{\F_q}},{\cal Q})$ (this cohomology group has dimension $1$ by Theorem \ref{HD}). Then we have 
\[\sum_{x\in U(\F_q)}t_{{\cal Q},q}(x)=(-1)^d\tau.
\]
In particular, Proposition \ref{FEV} holds true in this case. 
\end{enumerate}
\end{pr}
\proof{
$1.$ Let $X$ be an $\F_q$-scheme of finite type. For a morphism $f\colon Y\to X$ between $\F_q$-schemes of finite type and a $\Ql$-sheaf $\cal F$ on $X$, we have $t_{{\cal F},q}\circ f=t_{f^\ast{\cal F},q}$. For $\Ql$-sheaves ${\cal F},\cal G$ on $X$, we have $t_{{\cal F}\otimes{\cal G},q}=t_{{\cal F},q}\cdot t_{{\cal G},q}$ (see \cite[(1.1)]{Lau}). The assertion that $t_{{\cal Q},q}$ defines a quadratic form is verified by combining these general facts. 
It remains to  check that the map
\[t_{{\cal B}_{\cal Q},q}\colon U(\F_q)\times U(\F_q)\to
\Ql^\times
\]
is a perfect pairing. This can be proved by using Lang's isogeny. See \cite[Proposition A.18]{Bo} for the proof. 

$2.$ The displayed equality is a consequence of Grothendieck's trace formula and Theorem \ref{HD}. 
Proposition \ref{FEV} then follows from Proposition \ref{GS}, by noting that  $U$ is isomorphic to the perfection of $\mathbb{A}^d_{\F_q}$ as $\F_q$-schemes. 
\qed}

\vspace{5mm}

(Proof of Proposition \ref{FEV})

\vspace{3mm}

{We may replace $\F_q$ with a finite field extension. Therefore, by Corollary \ref{nondegdes}, we find $(V,{\cal Q}')$ with the following properties: $V$ is a perfect unipotent commutative $\F_q$-group scheme endowed with an abelian isogeny $f\colon U\to V$ and ${\cal Q}'$ is a non-degenerate quadratic $\Ql$-sheaf on $V$ such that $f^\ast{\cal Q}'\cong \cal Q$. 

Set $B:={\rm Aut}(U/V)$, which is naturally identified with $\mathop{\rm ker}\pi$. Then we have 
\[
f_\ast{\cal Q}\cong \bigoplus_{\psi\in B^\vee}{\cal Q}'\otimes{\cal L}_\psi
\]
where ${\cal L}_\psi$ denotes the multiplicative $\Ql$-sheaf corresponding to $\psi$. Note that ${\cal Q}'\otimes{\cal L}_\psi$ is a non-degenerate quadratic sheaf for any $\psi$. Then the assertion follows from Proposition \ref{nondegcase} by considering  the decomposition 
\[
H^d_c(U_{\overline{\F_q}},{\cal Q})\cong 
\bigoplus_{\psi\in B^\vee}H^d_c(V_{\overline{\F_q}},{\cal Q}'\otimes{\cal L}_\psi). 
\] 
\qed}

\subsection{Application: affine supersingular varieties}
\label{AffSS}

Let $q$ be a power of $p$ and $\F_q$ be a finite field with $q$ elements. Let us say that an affine smooth $\F_q$-variety $X$ is {\it supersingular} if any eigenvalue of ${\rm Fr}_q$ on $H^i_c(X_{\overline{\F_q}},\Ql)$ is of the form $\zeta\cdot\sqrt{q^i}$ for a root of unity $\zeta$. 
In this subsection, we introduce a certain class of affine supersingular varieties over $\F_q$  as an application of results established in the previous subsections. 

 Consider a pair $(C_0,U)$ with the following properties: 
\begin{itemize} 
\item $C_0$ is an affine connected smooth variety over $\F_q$. 
\item $U$ is a perfect unipotent commutative group scheme over $\F_q$. 
\item The perfection $C$ of $C_0$ is equipped with an 
$\F_q$-morphism 
$C\to U$ that is an abelian \'etale covering such that, for any character ${\rm Aut}(C/U)\to\Ql^\times$, the associated invertible $\Ql$-sheaf on $U$ is quadratic. 
\end{itemize}

\begin{pr}
Let $(C_0,U)$ be as above. Let $d=\dim U$. Then the  cohomology group $H^i_c(C_{0,\overline{\F_q}},\Ql)$ has  the following properties. 
\begin{enumerate}
\item It vanishes if $i<d$. 
\item If $i\geq d$, then every eigenvalue of the geometric Frobenius ${\rm Fr}_q$ on $H^i_c(C_{0,\overline{\F_q}},\Ql)$ is of the form $\zeta\cdot\sqrt{q^i}$ for a root of unity $\zeta$. 
\end{enumerate}
\end{pr}
\proof{
$1.$ This is a consequence of the affine Lefschetz theorem and the Poincar\'e duality. 

$2.$ Let $\pi\colon C\to U$ denote the covering map. Let $A:={\rm Aut}(C/U)$. We have the following  decomposition 
\[ \pi_\ast\Ql\cong\bigoplus_{\psi\in A^\vee}{\cal Q}_\psi,
\]
where $A^\vee$ denotes the character group and ${\cal Q}_\psi$ is the quadratic sheaf associated with $\psi$. Then we have 
\[
R\Gamma_c(C_{0,\overline{\F_q}},\Ql)\cong \bigoplus_{\psi\in A^\vee}R\Gamma_c(U_{\overline{\F_q}},{\cal Q}_\psi). 
\]
Then the assertion follows from Proposition \ref{FEV} together with the computation given in the subsection \ref{genrank}. \qed}

In the appendix \ref{exampleC}, we compute several examples of $(C_0,U)$ in low dimensional cases.

\subsection{$H^d_c(U_{\bar{k}},{\cal Q})$ vs.~$H^d(U_{\bar{k}},{\cal Q})$}

Let $\cal Q$ be a quadratic $\Lambda$-sheaf on $U$.
 Let $d=\dim U$.  
In this subsection, we prove that the canonical map $
H^d_c(U_{\bar{k}},{\cal Q})\to H^d(U_{\bar{k}},{\cal Q})
$ by forgetting the support is an isomorphism provided that $\cal Q$ is isogeneous. 
We do this by applying a result of Fourier--Deligne transform for unipotent groups, which we briefly recall below. 

Let $U$ be a perfect unipotent commutative group $k$-scheme of dimension $d$ and $\cal B$ be a bimultiplicative $\Lambda$-sheaf on $U\times_k U$ that is symmetric and non-degenerate. Consider the diagram 
\[
U\xleftarrow{{\rm pr}_1}U\times_kU\xrightarrow{{\rm pr}_2}U.
\]
Let $D_{\rm ctf}(U,\Lambda)$ denote the bounded derived category of constructible \'etale $\Lambda$-sheaves with finite tor-amplitude on $U$. We define the Fourier--Deligne transform 
\[\mathcal{F}_{\cal B}\colon D_{\rm ctf}(U,\Lambda)\to D_{\rm ctf}(U,\Lambda)
\]
by setting ${\cal F}_{\cal B}(K):=R{\rm pr}_{2!}({\rm pr}_1^\ast K\otimes^L{\cal B})[d]$. On the other hand, we define ${\cal F}_{{\cal B}\ast}(K)$ by setting $R{\rm pr}_{2\ast}({\rm pr}_1^\ast K\otimes^L{\cal B})[d]$. We have an obvious map ${\cal F}_{\cal B}(K)\to{\cal F}_{{\cal B}\ast}(K)$ of functors by forgetting the support. It is known that this map is in fact an isomorphism: 
\begin{thm}\label{!ast}(\cite[TH\'EOR\`EME 3.1]{Sai}, \cite[Theorem G.6]{BD})
The natural transformation ${\cal F}_{\cal B}(K)\to{\cal F}_{{\cal B}\ast}(K)$ constructed above is an isomorphism for any $K$. 
\end{thm}
\proof{
In \cite[TH\'EOR\`EME 3.1]{Sai}, this is proved by reducing to the case where $U=\G_{a,k}^{\rm perf}$, for which case it is proved in \cite{Lau}. On the other hand, in \cite[Theorem G.6]{BD}, the authors explain that the assertion can be deduced from the Fourier inversion formula. 
\qed}
\begin{pr}\label{HcH}
Let $U$ be a perfect unipotent commutative group $k$-scheme of dimension $d$. Let $\cal Q$ be an isogeneous quadratic $\Lambda$-sheaf on it. Then the canonical map 
\[ H_c^d(U_{\bar{k}},{\cal Q})\to H^d(U_{\bar{k}},{\cal Q})
\] 
is an isomorphism. 
\end{pr}
\proof{
We may assume that $k$ is algebraically closed. Then, by Corollary \ref{nondegdes}, we may find an \'etale isogeny $f\colon U\to V$ to a perfect unipotent commutative group $k$-scheme and a non-degenerate quadratic $\Lambda$-sheaf ${\cal Q}'$ on $V$ such that the pullback $f^\ast{\cal Q}'$ is isomorphic to $\cal Q$. Let $B$ be the automorphism group ${\rm Aut}(U/V)$ and $B^\vee$ be its character group. For $\psi\in B^\vee$, let ${\cal L}_\psi$ denote the multiplicative $\Lambda$-sheaf on $V$ corresponding to $\psi$.  We have a commutative diagram 
\[
\xymatrix{
H^d_c(U,{\cal Q})\ar[d]\ar[r]&H^d(U,{\cal Q})\ar[d]\\
\text{$\underset{\psi\in B^\vee}{\bigoplus}$}H_c^d(V,{\cal Q}'\otimes{\cal L}_\psi)\ar[r]&\text{$\underset{\psi\in B^\vee}{\bigoplus}$}H^d(V,{\cal Q}'\otimes{\cal L}_\psi). 
}
\]
Here the vertical arrows are isomorphisms induced  from $f_\ast{\cal Q}\cong \bigoplus_{\psi\in B^\vee}
{\cal Q}'\otimes{\cal L}_\psi$. Therefore, the assertion follows if we prove that each $H_c^d(V,{\cal Q}'\otimes{\cal L}_\psi)\to H^d(V,{\cal Q}'\otimes{\cal L}_\psi)$ is an isomorphism. In this way, we may reduce the assertion to the case where $\cal Q$ is non-degenerate. 

Suppose that $\cal Q$ is non-degenerate. Applying Theorem \ref{!ast} to ${\cal B}={\cal B}_{\cal Q}$ and $K=\cal Q$, we know that the map
\begin{equation}\label{comp!ast}
R{\rm pr}_{2!}({\rm pr}_1^\ast{\cal Q}\otimes{\cal B}_{\cal Q})\to R{\rm pr}_{2\ast}({\rm pr}_1^\ast{\cal Q}\otimes{\cal B}_{\cal Q})
\end{equation}
is an isomorphism. By the definition of ${\cal B}_{\cal Q}$, ${\rm pr}_1^\ast{\cal Q}\otimes{\cal B}_{\cal Q}$ is isomorphic to ${\rm pr}_2^\ast{\cal Q}^{-1}\otimes m^\ast{\cal Q}$, where $m\colon U\times_k U\to U$ is the map given by $(x,y)\mapsto x+y$. Since ${\cal Q}^{-1}$ is locally constant, the isomorphism (\ref{comp!ast}) can be identified with 
\[{\cal Q}^{-1}\otimes R{\rm pr}_{2!}(m^\ast{\cal Q})\to 
{\cal Q}^{-1}\otimes R{\rm pr}_{2\ast}(m^\ast{\cal Q}). 
\]
Thus the map $R{\rm pr}_{2!}(m^\ast{\cal Q})\to R{\rm pr}_{2\ast}(m^\ast{\cal Q})$ is an isomorphism. Since  the diagram
\[
\xymatrix{
U\times_kU\ar[r]^-{m}\ar[d]_-{{\rm pr}_2}&U\ar[d]\\
U\ar[r]&{\rm Spec}(k)
}
\]
is cartesian, the proper and generic base change theorems imply that the map $R\Gamma_c(U,{\cal Q})\to R\Gamma(U,{\cal Q})$ is an isomorphism, as desired. 
\qed}

\begin{appendices}

\section{Quadratic forms over finite abelian groups and their Gauss sums}

We review the theory of quadratic Gauss sums in a highly general form. In this appendix, 
for a finite abelian group $M$, we write $M^\vee$ for the character group ${\rm Hom}(M,\mathbb{C}^\times)$. 
\begin{dfA}\label{quadcl}
Let $M$ be a finite abelian group. 
\begin{enumerate}
\item We say that a mapping $Q\colon M\to \mathbb{C}^\times$ is {\rm a quadratic form} if the mapping $B_Q\colon M\times M\to \mathbb{C}^\times$ defined by $(x,y)\mapsto Q(x+y)Q(x)^{-1}Q(y)^{-1}$ is bilinear. 
\item We say that a quadratic form $Q$ on $M$ is {\rm non-degenerate} if the mapping $M\to M^\vee,x\mapsto B_Q(x,-),$ is an isomorphism. 
\item Let $Q$ be a quadratic form on $M$. We define  {\rm its Gauss sum} by setting 
\[
\tau_Q:=\sum_{x\in M}Q(x). 
\] 
\end{enumerate}
\end{dfA}
 Our definition of quadratic form is more general than that given in the literature, in that 
 we do not require that $Q(nx)$ should be $Q(x)^{n^2}$ for $n\in\mathbb{Z}$. This generality allows us to treat quadratic forms more flexibly, which will be seen below.  In what follows, we repeatedly use the following fact: for a quadratic form $Q$ on $M$, $B_Q$ is trivial if and only if it is a character. 
   
 Using the above definition, we can state the following (seemingly well-known) result for Gauss sums. 
\begin{prA}\label{GS}
Let $(M,Q)$ be a pair of a finite abelian group $M$ and a quadratic form $Q$ on it. Suppose that $Q$ is non-degenerate. Then the ratio 
\begin{equation*}
\frac{\tau_Q}{\sqrt{|M|}}
\end{equation*}
is a root of unity. 
\end{prA}
\begin{exA}\label{exeasy}
\begin{enumerate}
\item Suppose that $p$ is an odd prime number. Let $\psi\colon \F_p\to\mathbb{C}^\times$ be a non-trivial character. Then the mapping $Q_\psi\colon \F_p\to\mathbb{C}^\times, x\mapsto \psi(x^2),$ is a non-degenerate quadratic form and its Gauss sum is $\sum_{x\in\F_p}\psi(x^2).$ 
\item The mapping $Q\colon \F_2\to\mathbb{C}^\times$ given by $Q(0)=1, Q(1)=\sqrt{-1}$ is a non-degenerate quadratic form and its Gauss sum is $1+\sqrt{-1}$. 
\end{enumerate}
\end{exA}
 
The last of this section is devoted to the proof 
of  Proposition \ref{GS}.  
We start with noting the following property of $Q$. 
 \begin{lmA}\label{killcl}
 Let $M$ be a finite abelian group and $Q$ a quadratic form on it. Let $v\in M$ and $n\geq1$ be an integer which kills $v$. Then $Q(v)$ is an $n^2$-th root of unity. 
 \end{lmA}
 \proof{
 First note that $Q(0)=1$, since $B_Q(0,0)=1$. 
 Inductively on  an integer $m\geq1$, we can write 
 \[
 Q(mv)=Q(v)^m\cdot\prod_{i=1}^{m-1}B_Q(iv,v). 
 \]
 This implies that $Q(v)^n$ is an $n$-th root of unity. 
 \qed}
 
 The following lemma is useful to prove the proposition. 
 
 \begin{lmA}\label{twistcl}
Let $(M,Q)$ be the pair as in Proposition \ref{GS}. Let $\psi\colon M\to\C^\times$ be a group homomorphism. 
\begin{enumerate}
\item The mapping $M\to\C^\times, x\mapsto Q(x)\psi(x),$ is a non-degenerate quadratic form. 
\item The proposition holds true for $Q$ if and only if it holds true for $Q\psi$. 
\end{enumerate}
\end{lmA}
{\proof
$1.$ This is a simple computation. 

$2.$ Since $B_Q$ is a perfect pairing, there exists $a\in M$ such that $\psi(x)=B_Q(a,x)$. Then we can compute 
\begin{align*}
\sum_{x\in M}Q(x)\psi(x)&=\sum_{x\in M}Q(x)B_Q(a,x)\\
&=\sum_{x\in M}Q(a)^{-1}Q(a+x)=Q(a)^{-1}\sum_{x\in M}Q(x). 
\end{align*}
Since $Q(a)$ is a root of unity by Lemma \ref{killcl}, the assertion follows. 
\qed}

 In the following lemma,  we treat the simplest case of Proposition \ref{GS}. 
\begin{lmA}\label{vscl}
Let $(M,Q)$ be the pair as in Proposition \ref{GS}. 
Further suppose that $|M|=p$ for a prime number $p$. Then the statement of the proposition holds true. 
\end{lmA}
{\proof
We identify $M=\F_p$. 

First we treat the case when $p$ is odd. Let $\zeta_p=B_Q(\frac{1}{2},1)$. This is 
a primitive $p$-th root of unity by the non-degeneracy of $B_Q$. Let $Q_1\colon \F_p\to\C^\times$ be the quadratic form given by $x\mapsto \zeta_p^{x^2}$. 
Then we have $B_Q=B_{Q_1}$. In other words, there exists a character $\psi$ such that $Q=Q_1\psi$. Then, by Lemma \ref{twistcl}, the assertion for $Q$ is reduced to that for $Q_1$, which is classical. 

Next we consider the case $p=2$. Let $Q_2\colon\F_2\to\C^\times$ be the quadratic form given in Example \ref{exeasy}. It is straightforward to check that it is non-degenerate. Then we necessarily have $B_Q=B_{Q_2}$ since there is only one perfect pairing on $\F_2\times \F_2$ valued in $\C^\times$. 
Consequently, we can write $Q=Q_2\psi$ for a character $\psi$. Lemma \ref{twistcl} reduces the assertion for $Q$ to that for $Q_2$, which is straightforward. 
\qed}

\vspace{3mm}

(Proof of Proposition \ref{GS})

We argue by induction on $|M|$. If $|M|=1$, then the assertion is clear. Suppose that $|M|>1$. Choose a non-zero element $x\in M$ which is killed by a prime number $p$. First suppose that $B_Q(x,x)\neq1$. Then we have a direct sum decomposition 
(where $(-)^\perp$ is taken with respect to $B_Q$)
\begin{equation*}
M=\langle x\rangle\oplus\langle x\rangle^{\perp}. 
\end{equation*}
Since the restrictions of $Q$ to $\langle x\rangle,\langle x\rangle^\perp$ are non-degenerate, the assertion follows from the induction hypothesis and Lemma \ref{vscl}. 

Next suppose that $B_Q(x,x)=1$. This condition is equivalent to saying that $Q\colon \langle x\rangle\to\C^\times$ is a character. Let $\psi\colon M\to\C^\times$ be a character which extends $Q|_{\langle x\rangle}$. By Lemma \ref{twistcl}, we may replace $Q$ by $Q\psi^{-1}$ and may assume that $Q|_{\langle x\rangle}=1$.  

Suppose that $Q|_{\langle x\rangle}=1$. Then we can compute 
\begin{align*}
\tau_Q=\sum_{u\in M/\langle x\rangle}\sum_{v\in\langle x\rangle}Q(u+v)=\sum_{u\in M/\langle x\rangle}Q(u)\sum_{v\in\langle x\rangle}B_Q(u,v). 
\end{align*}
The sum $\sum_{v\in\langle x\rangle}B_Q(u,v)$ is zero unless $u\in \langle x\rangle^\perp/\langle x\rangle$. 
We have 
\[
\tau_Q=|\langle x\rangle|\sum_{u\in \langle x\rangle^\perp/\langle x\rangle}Q(u). 
\]
Note that the mapping $\langle x\rangle^\perp\to\C^\times, u\mapsto Q(u),$ induces a mapping 
$\langle x\rangle^\perp
/\langle x\rangle\to\C^\times$, which  defines a non-degenerate quadratic form on $\langle x\rangle^\perp
/\langle x\rangle$. Then the assertion follows from the induction hypothesis. 
\qed

When $M$ is an $\F_2$-vector space, the square $\tau_Q^2$ can be computed as follows. 
\begin{prA}
Let $M$ be a finite dimensional $\F_2$-vector space and $Q$ be a non-degenerate quadratic form on it. Let $B_Q$ be the associated pairing. 
\begin{enumerate}
\item The mapping $M\to \C^\times$ defined by $v\mapsto B_Q(v,v)$ is a character. 
\item There exists a unique element $a\in M$ which satisfies $B_Q(v,v)=B_Q(v,a)$ for any $v\in M$. The square $\tau_Q^2$ is equal to $Q(a)|M|$. 
\end{enumerate}
\end{prA}
{\proof 
$1.$ It follows since $B_Q(u,v)^2=B_Q(2u,v)=1$. 

$2.$ 
By the non-degeneracy of $B_Q$, there exists a unique element $a\in M$ such that $B_Q(u,a)=B_Q(u,u)$. 
We compute $\tau_Q^2$ as follows. 
\begin{align*}
\tau_Q^2=\sum_{u,v\in M}Q(u)Q(v)&=\sum_{u,v\in M}Q(u+v)B_Q(u,v)^{-1}\\
&=\sum_{w\in M}Q(w)\sum_{u\in M}B_Q(u,w-u)^{-1}\hspace{4mm}(w=u+v). 
\end{align*}
By replacing $B_Q(u,u)$ with $B_Q(u,a)$,  we can write 
\begin{align*}
\tau_Q^2=\sum_{w\in M}Q(w)\sum_{u\in M}B_Q(u,w-a)^{-1}. 
\end{align*}
The sum $\sum_u$ is zero unless $w=a$ and it is $|M|$ if $w=a$. Thus we get $\tau_Q^2=Q(a)|M|$. 
\qed}

\section{Low dimensional examples of affine supersingular varieties}\label{exampleC}

In this appendix, we give some examples of affine  supersingular varieties (in the sense of \ref{AffSS}) that appear in the context of quadratic sheaves. 
Let $k$ be a finite field of characteristic $p$; the constructions in what follows are valid for any perfect field, but we stick to the finite case as we define the supersingularity for affine varieties only in this case. 

\subsection{One dimensional examples}
Let $f(X),g_1(X),g_2(X)\in k[X]$ be  additive polynomials and let $a\in k$. We assume that the coefficient of $X$ in $f(X)$ is non-zero 
 and that $g_1(X),g_2(X)$ are non-zero. 
For such a datum, define the $k$-scheme $C$ by setting 
\[ C:={\rm Spec}(k[x,y]/(f(y)-g_1(x)g_2(x)-ax)). 
\]
\begin{prA}
The variety $C$ is a smooth supersingular curve. 
\end{prA}
\proof{
By the assumption on $f$, the derivation $d(f(y)-g_1(x)g_2(x)-ax)$ is of the form 
\[b\cdot dy+(\text{a term in which only $x$ appears})
\] for a non-zero $b\in k$. Then the smoothness follows from the Jacobian criterion.  To show the supersingularity, we may replace $k$ with a finite extension so that  all the roots of $f(X)$ are contained in $k$. Set $V:=\{y\in k\mid f(y)=0\}$.   Consider the map $\pi\colon C\to\G_{a,k}$ given by $(x,y)\mapsto x$. This map fits into the following cartesian diagram
\[
\xymatrix{
C\ar[r]^-\pi\ar[d]&\G_{a,k}\ar[d]^-h\\
\G_{a,k}\ar[r]^-f&\G_{a,k}
}
\]
where the map $h$ is given by $x\mapsto g_1(x)g_2(x)+ax$. 
Since all the roots of $f(Y)$ are in $k$, 
we have a decomposition $f_\ast\Ql\cong\bigoplus_{\psi\in V^\vee}{\cal L}_\psi$ over $\G_{a,k}$, where ${\cal L}_\psi$ is the $\psi$-isotipic part. Moreover, ${\cal L}_\psi$ are multiplicative as $f$ is a group homomorphism. Thus we have 
\[\pi_\ast\Ql\cong \bigoplus_{\psi \in V^\vee}h^\ast{\cal L}_\psi. \]
Note that the morphism $\G_{a,k}\times\G_{a,k}\to\G_{a,k}$ given by $(x,y)\mapsto h(x+y)-h(x)-h(y)$ is biadditive. From this observation, we know that 
$h^\ast{\cal L}_\psi$ are quadratic, hence the assertion. 
\qed}

When $f(Y)=Y^q-Y$ for some power $q$ of $p$, $g_1(X)=X$, and $a=0$, the curves are known as  the 
Van der Geer--Van der Vlugt curves (see \cite{GV}, \cite{BHetal}, and \cite{TT} for a detailed discussion on these curves). 

\subsection{Two dimensional examples}

In this subsection, we compute supersingular surfaces that derive from quadratic sheaves on $W_2$, the Witt scheme of length $2$. 

We start with reviewing some results on $W_2$. 
In the sequel, introducing the following notation is useful: 
\[\gamma(X_0,\dots,X_n):=\frac{\sum_{i=0}^nX_i^p-(\sum_{i=0}^nX_i)^p}{p}\in\mathbb{Z}[X_0,\dots,X_n]. 
\]
Recall that the ring structure of $W_2=\mathbb{A}^2_k$ is described as follows: 
\begin{align*}
&(x,y)+(z,w)=(x+z,y+w+\gamma(x,z)), \\
&(x,y)\cdot(z,w)=(xz,x^pw+z^py). 
\end{align*}
Using this, one can check that 
\[\sum_{i=0}^n(x_i,0)=(\sum_{i=0}^nx_i,\gamma(x_0,\dots,x_n) )
\]
by induction on $n$. As usual, we write $V\colon\G_{a,k}\to W_2$ for the Verschiebung and $F,R\colon W_2\to\G_{a,k}$ for the Frobenius and the restriction respectively.

The endomorphisms of $W_2$ can be classified as follows. 
\begin{lmA}\label{endW2}
A $k$-morphism $h\colon W_2\to W_2$ is additive if and only if $h(x,y)$ is of the form $(f(x),g_1(x)+g_2(y))$ for polynomials $f,g_1,$ and $g_2$ with the following properties: 
\begin{itemize}
\item $f(X)=\sum_{i=0}^nf_iX^{p^i}$ for some $f_i\in k$. 
\item $g_2(X)=f(X^{1/p})^p=\sum_{i=0}^nf_i^pX^{p^i}$. 
\item $g_1(X)=\gamma(f_0X^{p^0},f_1X^p,\dots,f_nX^{p^n})+R(X)$ where $R(X)\in k[X]$ is an additive polynomial. 
\end{itemize}
\end{lmA}
\proof{
We can write $h(x,y)=(\tilde{f}(x,y),g(x,y))$ for $\tilde{f},g\in k[x,y]$. 
Suppose that $h$ is additive. Then the composite map 
$W_2\xrightarrow{p}W_2\xrightarrow{h}W_2\xrightarrow{R}\G_{a,k}$ is zero. Since we have $p=VF$ and $F$ is faithfully flat, $R\circ h$ kills the group subscheme of $W_2$ defined by $V$. Thus 
there exists an additive polynomial $f(X)\in k[X]$ such that $\tilde{f}(x,y)=f(x)$. 

Comparing the second component of the equality $h((x,y)+(z,w))=h(x,y)+h(z,w)$, we get 
\begin{equation}\label{eqg}
g(x+z,y+w+\gamma(x,z))=g(x,y)+g(z,w)+\gamma(f(x),f(z)). 
\end{equation}
Putting $x,w=0$, we get $g(z,y)=g(0,y)+g(z,0).$ 
Therefore, one can write $g(x,y)=g_2(y)+g_1(x)$ for some polynomials $g_1,g_2$ with $g_1(0),g_2(0)=0$. 
Write $f(X)=\sum_{i=0}^nf_iX^{p^i}$. Then comparing the second component of $h(p(x,0))=ph(x,0)$ gives us the equality $g_2(X^p)=f(X)^p$. Consequently, $g_2$ is additive.

The equality (\ref{eqg}) is now reduced to the following 
\[g_1(x+z)-g_1(x)-g_1(z)=\gamma(f(x),f(z))-g_2(\gamma(x,z)).\]
On the other hand, the lemma below implies that the polynomial 
\[\widetilde{g}_1(X):=\gamma(f_0X^{p^0},f_1X^p,\dots,f_nX^{p^n})
\]
 satisfies the same relation as $g_1$. Thus $g_1(X)-\widetilde{g}_1(X)$ is additive, as claimed. 

Tracing back the discussion above, we can prove the only if part. 
\qed}
\begin{lmA}
Let $\widetilde{g}_1(X)$ be the element $\gamma(f_0X^{p^0}, f_1X^p,\dots,f_nX^{p^n})$ in 
the polynomial ring $\Z[X,f_0,\dots,f_n]$. Let $f(X):=\sum_{i=0}^nf_iX^{p^i}$. Then we have 
\[\widetilde{g}_1(X+Z)-\widetilde{g}_1(X)-\widetilde{g}_1(Z)\equiv\gamma(f(X),f(Z))-\sum_{i=0}^nf_i^p\gamma(X,Z)^{p^i}\hspace{3mm}{\rm mod}\ p. 
\]
\end{lmA}
\proof{
We can write 
\[\widetilde{g}_1(X+Z)-\widetilde{g}_1(X)-\widetilde{g}_1(Z)=\gamma(f(X),f(Z))-\sum_{i=0}^nf_i^p\cdot\frac{X^{p^{i+1}}+Z^{p^{i+1}}-(X+Z)^{p^{i+1}}}{p}\]
in $\Z[X,Z,f_0,\dots,f_n].$ Taking the $p$-th power of the congruence $(X+Z)^{p^i}\equiv X^{p^i}+Z^{p^i}\hspace{3mm}{\rm mod}\ p$, we know 
\[(X+Z)^{p^{i+1}}\equiv (X^{p^i}+Z^{p^i})^p\hspace{3mm}{\rm mod}\ p^2. 
\]
Then the assertion follows as $\gamma(X^{p^i},Z^{p^i})\equiv \gamma(X,Z)^{p^i}\hspace{3mm}{\rm mod}\ p. $
\qed
}

The computation in Lemma \ref{endW2} produces the following example of supersingular surfaces. 
\begin{prA}
Let $f_0,\dots,f_n$ be a sequence of elements in $k$ and $R(X)$ be an additive polynomial in $k[X]$. Define the polynomials $f(X),g_1(X),g_2(X)$ by the formulae listed in Lemma \ref{endW2}. 
\begin{enumerate}
\item Suppose that $p\neq2$. Then the $k$-scheme 
\[S:={\rm Spec}(k[x,y,z,w]/(x^p-x-zf(z),y^p-y+\gamma(x^p,-x)-x^pg_1(x)-x^pg_2(y)-f(x)^py))
\]
is an affine smooth supersingular surface. 
\item Suppose that $p=2$. Then the $k$-scheme 
\[S:={\rm Spec}(k[x,y,z,w]/(x^2-x-zf(z),y^2-y-x^2+x^3-x^2g_1(x)-x^2g_2(y)-f(x)^2y))
\]
is an affine smooth supersingular surface. 
\end{enumerate}
\end{prA}
\proof{
By Lemma \ref{endW2}, the map $h\colon W_2\to W_2$ defined by $h(x,y)=(f(x),g_1(x)+g_2(y))$ is additive. Let $L\colon W_2\to W_2$ be the map given by $L(x,y)=(x^p,y^p)-(x,y)$. This is an abelian isogeny with kernel canonically isomorphic to $\Z/p^2\Z$. Note that $L(x,y)$ is equal to $(x^p-x,y^p-y+\gamma(x^p,-x))$ if $p\neq2$ while it is $(x^2-x,y^2-y-x^2+x^3)$ if $p=2$. Then one can check that the scheme $S$ fits into the cartesian square 
\begin{equation*}
\xymatrix{
S\ar[d]\ar[r]^-\pi&W_2\ar[d]^-{\tilde{h}}\\
W_2\ar[r]^-{L}&W_2
}
\end{equation*}
where $\tilde{h}$ is given by $(x,y)\mapsto(x,y)\cdot h(x,y)$. Thus $S$ is a smooth surface. Moreover, we have 
\[
\pi_\ast\Ql\cong\bigoplus_{\psi\in(\Z/p\Z)^\vee}\tilde{h}^\ast{\cal L}_\psi, 
\]
where ${\cal L}_\psi$ is the locally constant sheaf on $W_2$ given by $\psi$. Since $\tilde{h}(a+b)-\tilde{h}(a)-\tilde{h}(b)$ is biadditive for $(a,b)\in W_2\times W_2$, we know that $\tilde{h}^\ast{\cal L}_\psi$ is quadratic. Then the assertion follows. 
\qed}

\section{Another group acting on $\cal Q$}\label{unitary}

Let $U$ be a perfect unipotent commutative group scheme over a perfect field $k$ and $\cal Q$ be a quadratic $\Lambda$-sheaf on $U$. In this appendix, we observe that one can attach to $\cal Q$ a certain group denoted by $Aut_U({\cal Q})$, which can be compared with the orthogonal groups of quadratic forms. 
\subsection{Definition of $Aut_U({\cal Q})$}

We define $Aut_U({\cal Q})$ as follows. 
Consider the presheaf 
\begin{equation*}
Aut_U\colon (Psch/k)^{\rm op}\to(\text{Sets})
\end{equation*}
defined by $Aut_U(S)=\{f\colon U_S\to U_S\mid\text{$f$ is a group automorphism}\}$. This is an \'etale sheaf of groups on $PSch/k$. 
We define $Aut_U({\cal Q})$ to be its sub-presheaf given by 
\[
Aut_U(Q)(S)=\{f\in Aut_U(S)\mid f^\ast{\cal Q}_S\cong {\cal Q}_S\}, 
\]
where ${\cal Q}_S$ deontes the pullback of $\cal Q$ to $U_S$; note that $f^\ast{\cal Q}_S$ is a quadratic sheaf since $f$ is a group homomorphism. The author does not know whether it is representable by a perfect group $k$-scheme. In this article, we can only prove the relative representability of the inclusion $Aut_U({\cal Q})\to Aut_U$. 
\begin{df}
A morphism ${\cal F}\to \cal G$ of presheaves on $PSch/k$ is {\rm relatively representable by a closed immersion} if, for any $S\in PSch/k$ and any $S\to{\cal G}$, the fiber product ${\cal F}\times_{\cal G}S$ is representable by a perfect $k$-scheme and the map 
${\cal F}\times_{\cal G}S\to S$ is a closed immersion. 
\end{df}

Using this notion, we state the result as follows. 
\begin{thm}\label{AutQ}
The injection $Aut_U({\cal Q})\to Aut_U$ is relatively representable by a closed immersion.  
\end{thm}

\begin{ex}
Let $V$ be a finite dimensional $k$-vector space and suppose that $U$ is the perfection of ${\rm Spec}({\rm Sym}^\bullet (V^\vee))$. Then there is a canonical injection of sheaves of groups $GL_V^{\rm perf}\to Aut_U$, where $(-)^{\rm perf}$ denotes the perfection functor. From this theorem, 
we get a perfect closed group subscheme $Aut_U({\cal Q})\times_{Aut_U}GL_V^{\rm perf}$ of $GL_V^{\rm perf}$. Note that the assignment $G\mapsto G^{\rm perf}$ gives a one to one correspondence between the smooth closed group subschemes of $GL_V$ and the perfect closed group subschemes of $GL_V^{\rm perf}$. Consequently, we know that there exists a smooth algebraic group $G$ over $k$ whose perfection is isomorphic to 
$Aut_U({\cal Q})\times_{Aut_U}GL_V^{\rm perf}$. 
\end{ex}

Let ${\cal B}_{\cal Q}$ be the bimultiplicative $\Lambda$-sheaf associated with $\cal Q$. To prove Theorem \ref{AutQ}, we also consider the presheaf $Aut_U({\cal B}_{\cal Q})$ given by 
\[Aut_U({\cal B}_{\cal Q})(S)=\{f\in Aut_V(S)\mid(f\times f)^\ast {\cal B}_{\cal Q}\cong {\cal B}_{\cal Q}\}. \]
\begin{lm}
The injection $Aut_U({\cal B}_{\cal Q})\to Aut_U$ is relatively representable by a closed immersion. 
\end{lm}
\proof{
Let $f\in Aut_U(S)$. We show that $Aut_U({\cal B}_{\cal Q})\times_{Aut_U}S$ is represented by a perfect closed subscheme. 
The bimultiplicative $\Lambda$-sheaf ${\cal B}_{\cal Q}$ corresponds to a group homomorphism $l_{\cal Q}\colon U_S\to U^\ast_S$ while $(f\times f)^\ast{\cal B}_{\cal Q}$ corresponds to $l_{f^\ast{\cal Q}}\colon U_S\to U^\ast_S$. The condition $(f\times f)^\ast{\cal B}_{\cal Q}\cong {\cal B}_{\cal Q}$ is equivalent to the condition that $l_{\cal Q}-l_{f^\ast{\cal Q}}=0$, i.e., $\mathop{\rm ker}(l_{\cal Q}-l_{f^\ast{\cal Q}})=U_S$. We note that $\mathop{\rm ker}(l_{\cal Q}-l_{f^\ast{\cal Q}})$ is  representable by a perfect closed group subscheme of $U_S$. Let $\cal U$ be the complement. Since $U$ is isomorphic to the perfection of a smooth $k$-scheme, the map $U\to{\rm Spec}(k)$ is universally open. Therefore, the image of $\cal U$ via $U_S\to S$ is open. Then the complement of this  image (with reduced subscheme structure) represents $Aut_U({\cal B}_{\cal Q})\times_{Aut_U}S$. 
\qed}

\vspace{5mm}

(Proof of Theorem \ref{AutQ})

\vspace{3mm}

{Since we already know that $Aut_U({\cal B}_{\cal Q})\to Aut_U$
 is relatively representable, it suffices to check that $Aut_U({\cal Q})\to Aut_U({\cal B}_{\cal Q})$ is relatively representable. Take $f\in Aut_U({\cal B}_{\cal Q})(S)$. Since we have $(f\times f)^\ast{\cal B}_{\cal Q}\cong {\cal B}_{\cal Q}$, the sheaf ${\cal M}:=f^\ast{\cal Q}\otimes
 {\cal Q}^{-1}$ is a multiplicative $\Lambda$-sheaf on $U_S$, which corresponds to a morphism $ S\to U^\ast$ of $k$-schemes. Then $S\times_{U^\ast}0$, where $0$ denotes the neutral  element of $U^\ast$,  represents $Aut_U({\cal Q})\times_{Aut_U({\cal B}_{\cal Q})}S$. 
 \qed}
 
 \subsection{Action of $Aut_U({\cal Q})$ on $C$}
 Let $\pi\colon C\to U$ be the cyclic covering that  trivializes $\cal Q$ and let $A:={\rm Aut}(C/U)$. 
 An element $f\in Aut_U({\cal Q})(k)$ acts on $C$ as follows.  Since we have $f^\ast {\cal Q}\cong {\cal Q}$, there exists an isomorphism $\alpha(f)'\colon C\to C$ such that the diagram
 \begin{equation*}
 \xymatrix{
 C\ar[r]^\pi\ar[d]_-{\alpha(f)'}&U\ar[d]^-f\\
 C\ar[r]_-\pi&U
 }
 \end{equation*}
 is commutative. If we moreover require that the induced map $\alpha(f)'\colon \pi^{-1}(0)\to\pi^{-1}(0)$ on the fibers of the neutral element $0$ of $U$ is the identity map, then there exists  a unique choice for $\alpha(f)'$. We write $\alpha(f)$ for the choice satisfying this normalization. From the normalization, one can check that $f\mapsto \alpha(f)$ defines an action of $Aut_U({\cal Q})(k)$ on $C$. 
 In the same way, $Aut_U({\cal Q})(\bar{k})$ acts on $C_{\bar{k}}$, where $\bar{k}$ is an algebraic closure of $k$. 
 
 From now on, we assume that $\cal Q$ is isogeneous and let $H$ denote the Galois group ${\rm Aut}(C_{\bar{k}}/U^\ast_{\bar{k}})$ of the covering $C_{\bar{k}}\to U^\ast_{\bar{k}}$ as in the subsection \ref{H}. We write $A$ for its center. 
 \begin{lm}\label{AutH}
 Let $\sigma\in H$ and $f\in Aut_U({\cal Q})(\bar{k})$. 
 \begin{enumerate}
 \item 
 The automorphism $\alpha(f)\sigma\alpha(f)^{-1}$ of $C_{\bar{k}}$ belongs to $H$. 
 \item If $\sigma$ is an element of the center $A$, then $\sigma$ and $\alpha(f)$ commute with each other. 
 \end{enumerate}
 \end{lm}
 \proof{
 $1.$ We have the following commutative diagram 
 \begin{equation*}
 \xymatrix{
 C_{\bar{k}}\ar[d]_-{\alpha(f)}\ar[r]^-\pi&U_{\bar{k}}\ar[d]_-f\ar[r]^-{l_{\cal Q}}&U^\ast_{\bar{k}}\\
 C_{\bar{k}}\ar[r]^-\pi&U_{\bar{k}}\ar[r]^-{l_{\cal Q}}&U^\ast_{\bar{k}}.\ar[u]^-{f^\ast}
 }
 \end{equation*}
 On the other hand, $\sigma$ fits into the commutative diagram 
 \begin{equation*}
 \xymatrix{
 C_{\bar{k}}\ar[d]_-\sigma\ar[r]&U^\ast_{\bar{k}}\ar[d]_-{\rm id}\\
 C_{\bar{k}}\ar[r]&U^\ast_{\bar{k}}. 
 }
 \end{equation*}
 The assertion follows from these commutative diagrams as we have $(f^{-1})^\ast=(f^\ast)^{-1}$. 
 
 $2.$ Let $\psi$ denote the injection $Z(H)\to\Lambda^\times$ corresponding to $\cal Q$. Then the map $Z(H)\to\Lambda^\times,\sigma\mapsto \psi(
 \alpha(f)\sigma\alpha(f)^{-1})$ corresponds to the locally constant sheaf $f^\ast{\cal Q}$. Then the condition $f^\ast{\cal Q}\cong{\cal Q}$ is equivalent to that $\psi(
 \alpha(f)\sigma\alpha(f)^{-1})=\psi(\sigma)$ for any $\sigma\in Z(H)$. As $\psi$ is injective, the assertion follows. 
 \qed}
 
 By this lemma, we can form the semidirect product 
 \[H\rtimes Aut_U({\cal Q})(\bar{k}), 
 \]
for which we write $G$ in the sequel. 
 This group acts on the both of $C_{\bar{k}}$ and $U_{\bar{k}}$ via $(\sigma,f)\mapsto\sigma\circ\alpha(f)$ and $(\sigma,f)\mapsto\sigma|_{U_{\bar{k}}}\circ f$ respectively. 
 
 Consider the direct sum decomposition 
 \[
 \pi_\ast\Lambda\cong\bigoplus_{\psi\in A^\vee}{\cal Q}_\psi, 
 \]
 where ${\cal Q}_\psi$ denotes the locally constant sheaf on $U_{\bar{k}}$ given by $\psi$.  Lemma \ref{AutH}.2 implies that the action of $G$ on $\pi_\ast\Lambda$ respects this decomposition. Consequently, $H^d_c(U_{\bar{k}},{\cal Q})$ gives a representation of the group $G$. 
 
 We end this appendix with observing that several geometric constructions of finite group representations in the literature can be regarded as special cases of the procedure explained above. 
 
  \begin{ex}\label{HWex}
 Let $q$ be a power of $p$.  
  For an \'etale sheaf $\cal L$ on $\G_{a,k}$ and a morphism of schemes $f\colon X\to\G_{a,k}$, we write ${\cal L}(f)$ for $f^\ast{\cal L}$. 
  \begin{enumerate}
  \item Let $\psi\colon \F_p\to\Ql^\times$ be an injective character and ${\cal L}_\psi$ be the  Artin--Schreier sheaf on $\G_{a,\F_q}$ as constructed in \ref{multGa}. Let $n$ be an integer $\geq1$. For 
  \[(x,y):=(x_1,\dots, x_n,y_1,\dots,y_n)\in\G_{a,\F_q}^{2n}, 
  \]
  we write $\langle x,y\rangle:=\sum_{i=1}^nx_iy_i$. Let ${\rm Fr}_q$ denote the $q$-th Frobenius map. Then 
  \[{\cal Q}:={\cal L}_\psi(\langle{\rm Fr}_q(x)-x,y\rangle)\]
   is an isogeneous quadratic $\Ql$-sheaf on $\G_{a,\F_q}^{2n}$. 
   The general linear group $GL_n(\F_q)$ is naturally regarded as a subgroup of $Aut_{\G_{a,\F_q}^{2n}}({\cal Q})(k)$ and $H^{2n}_c(\G_{a,\bar{\F_q}}^{2n},{\cal Q})$ gives a representation of $H\rtimes GL_{2n}(\F_q)$. This representation is studied in \cite{Tsu}. 
   \item Let $\F_{q,+}$ denote the group $\{z\in\F_{q^2}\mid z^q+z=0\}$, which we view as the Galois group of the covering $\G_{a,\F_{q^2}}\to\G_{a,\F_{q^2}}, z\mapsto z^q+z.$ Let $\psi\in \F_{q,+}^\vee$ be a non-trivial character and let ${\cal L}_\psi$ denote the corresponding locally constant sheaf on $\G_{a,\F_{q^2}}$. Then the sheaf 
   \[{\cal Q}:={\cal L}_\psi(\sum_{i=1}^nx_i^{q+1})
   \]
   on $\G_{a,\F_{q^2}}^n$ where $(x_1,\dots,x_n)$ is the standard coordinate functions, is an isogeneous quadratic $\Ql$-sheaf. The unitary group 
   \[U(n,+):=\{g\in GL_n(\F_{q^2})\mid \text{$g$ preserves the form $\sum_{i=1}^nx_i^{q+1}.$}\}
   \]
   is a subgroup of $Aut_{\G_{a,\F_{q^2}}^n}({\cal Q})(\F_{q^2})$. The representation  $H^n_c(\G_{a,\overline{\F_{q^2}}}^n,{\cal Q})$ of $H\rtimes U(n,+)$ is studied in \cite{IT}. 
    \end{enumerate}
  \end{ex}

\end{appendices}

\section*{Acknowledgement}
The author would like to thank Takahiro Tsushima for his interest in this work and many valuable comments. 
This project was started when the author learned Tsushima's result on the Suzuki curves at a conference which was organized by Prof.~Kenichi Bannai and his students. The author also thanks them for the organization. 
This work was supported  by 
RIKEN Special Postdoctoral Researcher Program.

\end{document}